\documentclass[german,a4paper,12pt]{article}
\usepackage{amsmath,amsthm,amsfonts,latexsym,amscd}

\usepackage{amsthm}

\usepackage{amssymb}
\swapnumbers
\theoremstyle{plain}

\theoremstyle{definition}

\theoremstyle{remark}
\begin{document}

\centerline{\large Differential forms and Hodge numbers}
\vspace{3mm}

\centerline{\large for toric complete intersections}
\vspace{3mm}

\centerline {\large Helmut A. HAMM}
\vspace{3mm}

{\bf Abstract:} We discuss conditions for complete intersections in a toric variety which allow to compute Hodge numbers if the complete intersection is a quasi-smooth complete variety. A preliminary step is the computation of the Euler characteristic of differential forms, we also look at symmetric or arbitrary forms instead of the usual alternating ones.\\

The notion of a complete intersection in a projective space is well-known. Here we will replace the projective space by an arbitrary toric variety and will use the notion of a complete intersection in this more general context. We will speak of toric complete intersections then (so a toric complete intersection need not be a toric variety). In contrast to the projective space a toric variety is not necessarily smooth. In fact, we will consider the case of a complete intersection which is smooth or at least quasi-smooth and therefore a $V$-manifold.\\

We will look at algebraic differential forms on such subvarieties, not only at the usual alternating forms but also at forms which are symmetric or arbitrary in the sense that no symmetry condition at all is imposed. \\

Since toric complete intersections are equipped with numerical data it is natural to compute the dimension of the cohomology groups of differential forms. We will restrict here to the easier task of computing the Euler characteristic.\\

In the case of alternating differential forms we will proceed to the question how to compute the Hodge numbers $h^{pq}$. The main ingredient is the computation of the Euler-Hodge characteristics $e^{pq}$ for non-degenerate complete intersections in tori. These invariants have already been computed by Danilov and Khovanski\u{\i} [D-K]; they dealt primarily with hypersurfaces and reduced the more general case of complete intersection to this special case by the trick of Lagrange multipliers. Here we will proceed in a different way: we will only partially reduce to hypersurfaces, in this way we will avoid increasing the dimension of the toric variety as in [D-K].\\

Using the Euler-Hodge characteristics for non-degenerate complete intersections in tori it is possible to compute the Hodge numbers $h^{pq}$ at least in two cases: for compact smooth toric complete intersections and for compact quasi-smooth varieties which can be decomposed into non-degenerate complete intersections in tori. We will introduce a class which includes both cases. To this end we intoduce the notion of non-degeneracy at infinity and show that the Hodge numbers can be computed for compact quasi-smooth varieties which admit a decomposition into smooth complete intersections which are non-degenerate at infinity. We conclude with several examples where we calculate the Hodge numbers.\\

Note that this paper constitutes a revised and enlarged version of the preprint ``Hodgezah\-len vollst\"andiger Durchschnitte in Tori" which was quoted in [H2].\\

{\bf 1. Toric varieties and subvarieties}\\

Note that for the computation of the Euler characteristic of differential forms (section 2) we will only need sections 1.1 - 1.3, whereas sections 1.4 - 1.6 prepare section 3.\\

{\bf 1.1.} A toric variety $X$ of dimension $m$ contains the $m$-dimensional torus $T\simeq (\mathbb{C}^*)^m$ as a dense open subset. Let $M\simeq\mathbb{Z}^m$ be the group of characters on $T$ and $N\simeq\mathbb{Z}^m$ the group of one-parameter subgroups of $T$. Then ${\cal O}(T)=\mathbb{C}[M]$ , because we prefer to work in the algebraic category. Note that there is a canonical pairing $<,>$ on $N\times M$ which comes from the canonical scalar product on $\mathbb{R}^m$. Now $X$ is defined by a fan $\cal F$ in $N_\mathbb{R}$, where $N_\mathbb{R}:=N\otimes_\mathbb{Z}\mathbb{R}\simeq \mathbb{R}^m$: $X=X_{\cal F}$. Note that $X$ is endowed with a $T$-action. See standard references like [K-K-M-S] or [O].\\

Remember that $X_{\cal F}$ is complete, smooth, quasi-smooth if and only if $\cal F$ is complete, regular, or simplicial, respectively.\\ 

Let us fix a fan $\cal F$ and put $X:=X_{\cal F}$. We may describe $X$ in a different way, namely as a quotient. This possibility has been remarked very early, see [Dz], but with restrictions on the fan; for the general case see [C]. See [A], [Ba-C], too. The model is $\mathbb{P}_m\simeq(\mathbb{C}^{m+1}\setminus\{0\})/\sim$ . In [Ba-C] this method has been used in order to study the Hodge theory of toric hypersurfaces. \\

We will use a modified approach here, however, which is comparable to the use of the graph of a mapping instead of the mapping itself and has numerical advantages. For the relation to the general process of forming quotients see also [H3].\\

So we proceed as follows (see also [H3]): Let $p_1,\ldots,p_r$ be the generators of the semigroups $\sigma\cap N$, $\sigma$ being the edges (i.e. one-dimensional cones) in $\cal F$. Now let us look at the following fan ${\cal F}'$ in $\mathbb{R}^{r+m}$ (and not in $\mathbb{R}^r$, as usual): Instead of $p_j$ let us take the canonical unit vector $e_j$. If $\sigma$ is a cone in $\cal F$, generated by $p_{j_1},\ldots,p_{j_s}$, we have a corresponding cone $\sigma'$, generated by $e_{j_1},\ldots,e_{j_s}$. Let us take all cones $\sigma'$ obtained in this way, together with their faces. In this way we get a fan ${\cal F}'$. Let $X':=X_{{\cal F}'}$.\\ 

We have a partition of $\mathbb{C}^{r+m}$ by subsets of the form $S_1\times\ldots\times S_{r+m}$ with $S_j=\mathbb{C}^*$ or $S_j=\{0\}$. If $\sigma'$ is a cone, generated by $e_{j_1},\ldots,e_{j_s}$, then $X_{\sigma'} =U_1\times\ldots\times U_{r+m}$ with $U_j := \mathbb{C}$ for $j=j_1,\ldots,j_s$ and $U_j:= \mathbb{C}^*$ otherwise. Now $X_{\sigma'}$ contains exactly one closed $(\mathbb{C}^*)^{r+m}$-orbit $O_{\sigma'}$, and $O_{\sigma'}=S_1\times\ldots\times S_{r+m}$ with $S_j=\{0\}$ f\"ur $j=j_1,\ldots,j_s$ and $S_j=\mathbb{C}^*$ otherwise, the partition above is therefore the partition into orbits. Furthermore, $X_{{\cal F}'} = \bigcup_{\tau\in {\cal F}'} X_\tau$ is an open subset of $\mathbb{C}^r\times(\mathbb{C}^*)^m$ which can be written as union of sets of the partition above. Therefore $X_{{\cal F}'}$ is much more intuitive than $X_{\cal F}$.\\

We have an action of the torus $(\mathbb{C}^*)^r$ on $\mathbb{C}^{r+m}$: $$c\circ(\zeta_1,\ldots,\zeta_r,z_1,\ldots,z_m):= (c_1\zeta_1,\ldots,c_r\zeta_r,c^{-p^1} z_1,\ldots,c^{-p^m}z_m),$$ 
where $p^i:=(p_{1i},\ldots,p_{ri})$. In particular, $(\mathbb{C}^*)^r$ acts on $X_{{\cal F}'}$, too.\\

Let $\sigma$ be the cone in $\cal F$ generated by $p_{j_1},\ldots,p_{j_s}$. Then we have a homomorphism from ${\cal O}(X_{\sigma}) = \mathbb{C}[\mathbb{Z}^m\cap\check{\sigma}]$ to
${\cal O}(X_{\sigma'}) = \mathbb{C}[\mathbb{Z}^{r+m}\cap\check{\sigma'}]$, defined by  $z^q\mapsto \zeta_1^{<p_1,q>}\cdot\ldots\cdot \zeta_r^{<p_r,q>} z^q$. For if  $q\in \mathbb{Z}^m$, $q\in \check{\sigma}$, i.e. $<p_j,q>\geq 0$ for $j=j_1,\ldots,j_s$, then $q$ corresponds to the element $z^q\in 
\mathbb{C}[\mathbb{Z}^m\cap\check{\sigma}]$. Obviously, $\zeta_1^{<p_1,q>}\cdot\ldots\cdot \zeta_r^{<p_r,q>} z^q\in \mathbb{C}[\mathbb{Z}^{r+m}\cap\check{\sigma'}]$. So we have a morphism $\pi_\sigma\colon X_{\sigma'}\rightarrow X_\sigma$. Since $(\mathbb{C}^*)^r$ is an abelian and therefore reductive group which acts on the affine variety $X_{\sigma'}$ the corresponding algebraic (and therefore categorical) quotient $X_{\sigma'}//(\mathbb{C}^*)^r := Spec ({\cal O}(X_{\sigma'})^{(\mathbb{C}^*)^r})$ exists. Obviously, the morphisms  $\pi_{\sigma}$ fit together to a morphism $\pi\colon X_{{\cal F}'}\rightarrow X_{\cal F}$.\\

{\bf Theorem 1.1.1·} (cf. [H3]) If $\sigma\in \cal F$, the morphism $\pi_{\sigma}$ induces an isomorphism of the algebraic quotient $X_{\sigma'}//(\mathbb{C}^*)^r$ onto $X_\sigma$.\\

{\bf Proof.} Let $\sigma$ be as above. Then the monomials $\zeta_1^{<p_1,q>}\cdot\ldots\cdot \zeta_r^{<p_r,q>} z^q$, $q\in \mathbb{Z}^m$, $<p_j,q>\geq 0$ f\"ur $j=j_1,\ldots,j_s$, form a basis of the vector space of the $(\mathbb{C}^*)^r$-invariant elements of ${\cal O}(X_{\sigma'})$.\\

Therefore $X_{\cal F}$ is a categorical quotient, it can be obtained by glueing (sc. affine) algebraic quotients. In fact, if $\cal F$ is simplicial we have a geometric quotient, i.e. the (closed) points of $X_{\cal F}$ correspond to the $(\mathbb{C}^*)^r$-orbits:\\

{\bf Theorem 1.1.2:} (cf. [H3]) If ${\cal F}$ is simplicial $X_{\cal F}$ is a geometric quotient with respect to the $(\mathbb{C}^*)^r$-action on $X_{{\cal F}'}$.\\

We can describe the relation between ${\cal O}_{X'}$ and ${\cal O}_X$ in the following way, too. The action of $(\mathbb{C}^*)^r$ defines a $\mathbb{Z}^r$-grading on $\pi_*{\cal O}_{X'}$: $\pi_*{\cal O}_{X'}=\oplus_{s\in\mathbb{Z}^r} (\pi_*{\cal O}_{X'})_s$. If $\sigma$ is a cone as in Theorem 1.1.1,  $(\pi_*{\cal O}_{X'})_s(X_{\sigma})$ is the complex vector space spanned by all monomials $\zeta_1^{s_1+<p_1,q>}\cdot\ldots\cdot \zeta_r^{s_r+<p_r,q>} z^q$ with $q\in \mathbb{Z}^m$, $<p_j,q>\geq -s_j$ f\"ur $j=j_1,\ldots,j_s$. Therefore we can speak of graded $\pi_*{\cal O}_{X'}$-modules ${\cal T}$: ${\cal T}=\oplus_{s\in\mathbb{Z}^r} {\cal T}_s$.\\

In particular, we obtain ${\cal O}_X=(\pi_*{\cal O}_{X'})_0$, see the proof of Theorem 1.1.1; so  $(\pi_*{\cal O}_{X'})_s$ is an ${\cal O}_X$-module. \\

{\bf 1.2.} The advantage of representing $X$ as a quotient becomes clear when we pass to subvarieties of $X$. We will do this first in a special case.\\

This is motivated by the case of projective varieties: A subvariety $V$ of $\mathbb{P}_m(\mathbb{C})$ can be written as $\{x\in\mathbb{P}_m(\mathbb{C})\,|\,{\rm f}_1(x)=\ldots={\rm f}_k(x)=0\}$, where ${\rm f}_j$ is a section of the line bundle
${\cal O}_{\mathbb{P}_m(\mathbb{C})}(d_j)$ on $\mathbb{P}_m(\mathbb{C})$. We can avoid the language of line bundles: ${\rm f}_j$ corresponds to a homogeneous polynomial $F_j$ on $\mathbb{C}^{m+1}$; then $V=\pi(V')$, where $V':=\{z\in\mathbb{C}^{m+1}\setminus\{0\}\,|\,F_1(z)=\ldots=F_k(z)=0\}$, and $\pi:\mathbb{C}^{m+1}\setminus\{0\}\longrightarrow \mathbb{P}_m(\mathbb{C})$ is the canonical mapping.\\  

For $s\in\mathbb{Z}^r$ we have the following action of $(\mathbb{C}^*)^r$ on $\mathbb{C}$: $c\circ t:=c_1^{s_1}\cdot\ldots\cdot c_r^{s_r}t$. Remember that $(\mathbb{C}^*)^r$ acts on $\mathbb{C}^{r+m}$ as indicated before.\hfill(*)\\

Let $H\in \mathbb{C}[\zeta_1,\ldots,\zeta_r,z_1,z_1^{-1},\ldots,z_m,z_m^{-1}]$. Obviously:\\

{\bf Lemma 1.2.1:} The following conditions are equivalent:\\
a) $H$ is equivariant with respect to (*),\\
b) $H$ is a linear combination of monomials of the form $\zeta^\rho z^q$ with $\rho_j-<p_j,q>=s_j$, $j=1,\ldots,r$,\\
c) $H$ defines a global section $\rm h$ of $(\pi_*{\cal O}_{X'})_s$.\\

 We can also proceed in the following way. Let us consider  the following subsheaf ${\cal S}_s$ of $\iota_*{\cal O}_T$, where $\iota\colon T\rightarrow X_{\cal F}$ denotes the inclusion of $T$ in $X_{\cal F}$: if $\sigma\in {\cal F}$ is generated by $p_{j_1},\ldots,p_{j_l}$, ${\cal S}_s(X_\sigma)$ is by definition generated by the $z^q$ with $<p_{j_\lambda},q>\ge -s_{j_\lambda}, \lambda=1,\ldots,l$. Here $z^q\in{\cal O}(T)$ is simultaneously considered as a section in $(\iota_*{\cal O}_T)(X_\sigma)$.\\

 If $H$ satisfies the conditions of Lemma 1.2.1 and if we replace $\zeta^\rho z^q$ by $z^q$ within $H$, we obtain a Laurent polynomial $h\in \mathbb{C}[z_1,z_1^{-1},\ldots,z_m,z_m^{-1}]$, and $h$ can be viewed as a section in ${\cal S}_s$. On the other hand, each Laurent polynomial which is a linear combination of monomials $z^q$ with $<p_j,q>\ge -s_j$, $j=1,\ldots,r$, can be interpreted as a section of ${\cal S}_s$ and leads to an equivariant Laurent polynomial $H$, replacing $z^q$ by $\zeta^\rho z^q$ with $\rho_j:=<p_j,q>+s_j$, $j=1,\ldots,r$.\\

 In particular the global sections of $(\pi_*{\cal O}_{X'})_s$ and ${\cal S}_s$ correspond to each other. It is easy to see that in fact $(\pi_*{\cal O}_{X'})_s\simeq{\cal S}_s$.\\

For $i=1,\ldots,k$  fix $d_{i1},\ldots,d_{ir}\in\mathbb{Z}$. We can look at the following action of $(\mathbb{C}^*)^r$ on $\mathbb{C}^k$: $c\circ (t_1,\ldots,t_k):= (c_1^{d_{11}}\cdot\ldots\cdot c_r^{d_{1r}}t_1,\ldots,c_1^{d_{k1}}\cdot\ldots\cdot c_r^{d_{kr}}t_k)$. Remember that $(\mathbb{C}^*)^r$ acts on $\mathbb{C}^{r+m}$ as indicated before.\\

 Let $G_1,\ldots,G_k\in \mathbb{C}[\zeta_1,\ldots,\zeta_r,z_1,z_1^{-1},\ldots,z_m,z_m^{-1}]$, $G:=(G_1,\ldots,G_k)$. Suppose that $G$ is equivariant. Then the set $Y':=\{(\zeta,z)\in X'\,|\,G(\zeta,z)=0\}$ is $(\mathbb{C}^*)^r$-invariant. Define $Y:=\pi(Y')$.\\

We want to show that $Y$ is a closed algebraic subspace of $X$. This is easy under a supplementary hypothesis:\\

 Suppose now that the following holds: For every $\sigma\in{\cal F}$ generated by $p_{j_1},\ldots,p_{j_l}$ there is a $q\in M$ such that $<p_{j_\lambda},q>=-s_{j_\lambda}, \lambda=1,\ldots,l$. Then $z^q\in{\cal S}_s(X_\sigma)$, and $z^q$ trivializes ${\cal S}_s|X_\sigma$: ${\cal S}_s$ is a line bundle. In particular we may then speak of zeroes of sections.\\

 Put $d_i:=(d_{i1},\ldots,d_{ir}), i=1,\ldots,k$. Then the procedure above associates to $G_1,\ldots,G_k$ Laurent polynomials $g_1,\ldots,g_k\in \mathbb{C}[z_1,z_1^{-1},\ldots,z_m,z_m^{-1}]$ as well as global sections ${\rm g}_1$ in $(\pi_*{\cal O}_{X'})_{d_1}$, $\ldots,{\rm g}_k$ in $(\pi_*{\cal O}_{X'})_{d_k}$. Obviously $g_i$ is a linear combination of monomials $z^q$ with $<p_j,q>\ge-d_{ij}$, $j=1,\ldots,r$.\\

{\bf Lemma 1.2.2:} Assume that the following holds: For every $\sigma\in{\cal F}$, spanned by $p_{j_1},\ldots,p_{j_s}$, and every $i\in\{1,\ldots,k\}$  there is  a $q\in M$ such that $<p_{j_\lambda},q>= -d_{ij_\lambda}$, $\lambda=1,\ldots,s$. Then $Y$ is the set of zeroes of the sections ${\rm g}_1,\ldots,{\rm g}_k$,\\ 

Note that the hypothesis yields that we have sections in line bundles, so we can speak of zeroes of these sections. In particular, $Y$ is a closed algebraic subset of $X$.\\

{\bf Lemma 1.2.3:} Assume that $\cal F$ is simplicial or that the hypothesis of Lemma 1.2.2 holds. Then $Y'=\pi^{-1}(Y)$.\\

{\bf Proof:} If $\cal F$ is simplicial the statement is obvious because $Y'$ is $(\mathbb{C}^*)^r$-invariant and we have a geometric quotient.\\
Now assume the hypothesis of Lemma 1.2.2. Let us look at $\pi|\pi^{-1}(X_\sigma):\pi^{-1}(X_\sigma)\to X_\sigma$. Choose $q$ as above, then $Y'$ is defined by the invariant polynomials $z^{-q}G_j$, $j=1,\ldots,k$. Therefore we have $Y'\cap \pi^{-1}(X_\sigma)=\pi^{-1}(\pi(Y'\cap\pi^{-1}(X_\sigma)))$.\\

 For the next section let us make the following remark concerning non-simplicial fans. The hypothesis of Lemma 1.2.3 implies the following condition on the numbers $d_{ij}$: \\

For every $\sigma\in{\cal F}$, spanned by $p_{j_1},\ldots,p_{j_s}$, and every $i\in\{1,\ldots,k\}$, the vector spaces spanned by $p_{j_1},\ldots,p_{j_s}$ resp. $(p_{j_1},d_{ij_1}),\ldots,(p_{j_s},d_{ij_s})$ have the same dimension.\hfill(**)\\

{\bf 1.3.} Now we drop any hypothesis about the fan.\\

{\bf Lemma 1.3.1:}\\
a) $Y$ is a closed algebraic subset of $X$.\\
b) $Y\cap T=\{z\in T\,|\,g_1(z)=\ldots=g_k(z)=0\}$.\\

{\bf Proof:} a) We know that $\pi|\pi^{-1}(X_\sigma):\pi^{-1}(X_\sigma)\to X_\sigma$ is an algebraic quotient. Since $Y'$ is a closed $(\mathbb{C}^*)^r$-invariant algebraic subset of $\pi^{-1}(X_\sigma)$ we conclude by [Kr] p. 96 that $Y\cap X_\sigma$ is a closed algebraic subset of $X_\sigma$.\\
b) obvious.\\

Let $Y'\neq\emptyset$. Generalizing the projective case let us call $Y$ a (toric) complete intersection if $Y$ has codimension $k$. We have:\\

{\bf Lemma 1.3.2:} a) If $Y$ has codimension $k$ in $X$ we have that $Y'$ has codimension $k$ in $X'$, too.\\
b) Under assumption (**) the converse holds.\\

{\bf Proof:} First we remark the following: Let $\sigma\in{\cal F}$ and $\tau\in{\cal F}'$ such that $\pi(O_\tau)=O_\sigma$. Then the fibres of the mapping $\pi|O_\tau:O_\tau\to O_\sigma$ have the dimension $\dim\,O_\tau-\dim\,O_\sigma=(m+r-\dim\,\tau)-(m-\dim\,\sigma)=r-\dim\,\tau+\dim\,\sigma$.

\vspace{2mm}
a) Let $\sigma, \tau$ be as above. Then $\dim\,\tau\ge\dim\,\sigma$, hence $\dim\,Y'\cap O_\tau\le (m-k)+(r-\dim\,\tau+\dim\,\sigma)\le m+r-k$.

\vspace{2mm}
b) First let us assume that we have the hypothesis of Lemma 1.2.3, so $Y'=\pi^{-1}(Y)$. Let $\sigma$ be a cone of ${\cal F}$. Assume that $\dim\,Y\cap O_\sigma>m-k$: Choose a face $\tau$ of $\sigma'$ such that $\dim\,\tau=\dim\,\sigma$, then $\pi(O_\tau)=O_\sigma$. Since $Y'=\pi^{-1}(Y)$ we have that $\dim\,Y'\cap O_\tau=\dim\,Y\cap O_\sigma+(r-\dim\,\tau+\dim\,\sigma)>m+r-k$, which gives a contradiction.\\
In general we obtain the hypothesis of Lemma 1.2.2 if we replace $M$ by some lattice $\tilde{M}$ which contains $M$ as a sublattice of finite index. This gives a toric variety $\tilde{X}$ with a finite map $p:\tilde{X}\to X$. In analogy to $\pi:X'\to X$ we can form $\tilde{\pi}:\tilde{X}'\to\tilde{X}$, and we have a finite mapping $p':\tilde{X}'\to X'$ such that $\pi\circ p'=p\circ\tilde{\pi}$. Let $\tilde{Y}'$ and $\tilde{Y}$ be defined as $Y'$ and $Y$ but starting from our different lattice. Since $p'$ is finite we have $\dim\tilde{Y}'\le m+r-k$. By the special case treated before we get $\dim\tilde{Y}\le m-k$. Since $p$ is finite we obtain $\dim\,Y\le m-k$, too.\\

{\bf Remark:} We cannot drop the assumption (**) in Lemma 1.3.2: Let $\sigma$ be a cone of $\cal F$ which is not simplicial, spanned by $p_{j_1},\ldots,p_{j_l}$. Then $Y':=\{\zeta_{j_1}=\ldots=\zeta_{j_l}=0\}$ defines a complete intersection in $X'$ of codimension $l$, where $d_i=e_{j_i}, i=1,\ldots,l$, whereas $Y=\overline{O_\sigma}$ has codimension $\dim\sigma<l$.\\

{\bf Lemma 1.3.3:} Assume that the Jacobian matrix of $G$ has rank $k$ along $Y'$.\\
a) If ${\cal F}$ is simplicial $Y$ is a $V$-manifold of pure dimension $m-k$.\\
b) If ${\cal F}$ is regular $Y$ is a manifold of pure dimension $m-k$.\\

{\bf Proof:} a) It is sufficient to show that $Y\cap X_\sigma$ is a purely $m-k$-dimensional $V$-manifold around $Y\cap O_\sigma$, where $\sigma\in{\cal F}$ is arbitrary. Without loss of generality we may assume that $\sigma$ is generated by $p_1,\ldots,p_s$ and that $p_1,\ldots,p_s,e_{s+1},\ldots,e_m$ is a basis of $\mathbb{R}^m$. Let $A:=\mathbb{C}^s\times\{(1,\ldots,1)\}\times(\mathbb{C}^*)^{m-s}\subset\mathbb{C}^r\times(\mathbb{C}^*)^m$ and $\Gamma$ the isotropy group of $A$.\\
 Then $(\mathbb{C}^*)^r\circ A=X_{\sigma'}$: Let $(\zeta,z)\in X_{\sigma'}$. Then we have to show the existence of a $c\in(\mathbb{C}^*)^r$ such that the following equations hold:
$$c_{s+1}\zeta_{s+1}=\ldots=c_r\zeta_r=1,
c_1^{-p_{1j}}\cdots c_r^{-p_{rj}}z_j=1,j=1,\ldots,s$$
Note that $c_{s+1},\ldots,c_r$ are determined by the first equations. Let $\log$ denote a fixed branch of the logarithm, e.g. the standard branch $\log\,re^{i\phi}:=\ln r+i\phi, -\pi<\phi\le\pi$, $\gamma_j:={\log c_j\over 2\pi i}, j=s+1,\ldots,r$. Then let $(\gamma_1,\ldots,\gamma_s)$ be a solution of the system of linear equations:
$$-p_{1j} \gamma_1-\ldots-p_{rj}\gamma_r +\log z_j=0 , j=1,\ldots,s$$ 
Now put $c_j:=e^{2\pi i\gamma_j}$, $j=1,\ldots,s$.\\
 The group $\Gamma$ is finite: 
Look at all $(\gamma_1,\ldots,\gamma_s)$ such that 
$$-p_{1j} \gamma_1-\ldots-p_{sj}\gamma_s\in\mathbb{Z}, j=1,\ldots,s$$
Let $D:=\det (p_{ij})_{1\le i,j\le s}$. Then $D\gamma_j\in\mathbb{Z}$ for all $j$ by Cramer's rule. So $G$ is a subgroup of $C_D^s\times\{(1,\ldots,1)\}$, where $C_D\subset \mathbb{C}^*$ is the cyclic group of order $D$.\\
 Since $Y\cap X_\sigma\simeq (A\cap Y')/\Gamma$ it is sufficient to show that $A\cap Y'$ is smooth along $A\cap Y'\cap O_{\sigma'}$. Now
$$\zeta_j{ \partial G_i\over\partial\zeta_j}-p_{j1}  z_1{\partial G_i\over\partial z_1}-\ldots-p_{jm}  z_m{\partial G_i\over\partial z_m}=d_{ij} G_i,\hfill{*}$$
$i=1,\ldots,k$, $j=1,\ldots,r$. Let $(\zeta,z)\in A\cap Y'\cap O_{\sigma'}$, in particular $\zeta_1=\ldots=\zeta_s=0$, $\zeta_{s+1}=\ldots=\zeta_r=1$. Then:
$$-p_{j1}  z_1{\partial G\over\partial z_1}-\ldots-p_{jm}  z_m{\partial G\over\partial z_m}=0,\,\,j=1,\ldots,s,$$
$${\partial G\over\partial\zeta_j}-p_{j1}  z_1{\partial G\over\partial z_1}-\ldots-p_{jm}  z_m{\partial G\over\partial z_m}=0,\,\,j=s+1,\ldots,r.$$
Since $\det((p_{ij})_{1\leq i,j\leq s})\neq 0$ the vectors ${\partial G\over\partial\zeta_{s+1}},\ldots,{\partial G\over\partial\zeta_r},{\partial G\over\partial z_1},\ldots,{\partial G\over\partial z_s}$ depend at $(\zeta,z)$ linearly on ${\partial G\over\partial z_{s+1}},\ldots,{\partial G\over\partial z_m}$. Consequently, by hypothesis the partial derivatives ${\partial G\over\partial\zeta_1},\ldots,{\partial G\over\partial\zeta_s},$ ${\partial G\over\partial z_{s+1}},\ldots,{\partial G\over\partial z_m}$ in $(\zeta,z)$ must span a space of dimension $k$, which implies our assertion, because $\zeta_1,\ldots,\zeta_s,z_{s+1},\ldots,z_m$ are coordinates for $A$.

\vspace{2mm}
 b) We can assume above: $p_1=e_1,\ldots,p_s=e_s$. Then $\Gamma$ is trivial.\\ 

 Under the assumption of Lemma 1.3.3a)  resp. b) we call $Y$ a quasi-smooth resp. smooth complete intersection.\\

 The hypothesis of Lemma 1.3.3 should not be confused with the stronger condition that the Jacobian matrix of $G|O_{\sigma'}$ has rank $k$ along $Y'\cap O_{\sigma'}$ for all $\sigma\in \cal F$:\\

 For $\sigma\in {\cal F}$, generated by $p_{j_1},\ldots,p_{j_s}$, let
\\
  $g_i^\sigma := \sum_{q\in M\,:\,<p_{j_1},q>=-d_{ij_1},\ldots, <p_{j_s},q>=-d_{ij_s}} a_{iq} z^q$.\\

{\bf Lemma 1.3.4.} The following conditions are equivalent:\\
a) the Jacobian matrix of $G|O_{\sigma'}$  has rank $k$ along $Y'\cap O_{\sigma'}$,\\
b) the mapping $(g_1^\sigma,\ldots,g_k^\sigma)\colon (\mathbb{C}^*)^m\rightarrow \mathbb{C}^k$ has no critical point $z$ with $g_1^\sigma(z)=\ldots=g_k^\sigma(z)=0$.\\

 {\bf Proof:} Without loss of generality we may assume that $\sigma$ is generated by $p_1,\ldots,p_s$. Then $O_{\sigma'}=\{0\}\times(\mathbb{C}^*)^{m+r-s}\subset (\mathbb{C}^*)^{r+m}$.\\
We can identify $/\mathbb{C}^*)^m$ with $A:=\{0\}\times\{(1,\ldots,1)\times(\mathbb{C}^*)^m\subset (\mathbb{C}^*)^s\times(\mathbb{C}^*)^{r-s}\times(\mathbb{C}^*)^m\simeq(\mathbb{C}^*)^{r+m}$. Under this identification, $g_i^\sigma$ corresponds to the restriction of $G_i$. We have equation (*) of the proof of Lemma 1.3.3 again. Therefore $\frac{\partial G}{\partial\zeta_{s+1}},\ldots,\frac{\partial G}{\partial\zeta_{r}}$ are linearly dependent on $\frac{\partial G}{\partial z_1},\ldots,\frac{\partial G}{\partial z_m}$. So we have that the following statements are equivalent:\\
(i) the Jacobian matrix of $G|O_{\sigma'}$  has rank $k$ along $A$,\\
(ii) the Jacobian matrix of $G|A$  has rank $k$ along $A$.\\
Now $O_{\sigma'}$ is the $(\mathbb{C}^*)^{r+m}$-orbit of $A$, so a) is equivalent to (i), whereas b) is equivalent to (ii).\\

{\bf 1.4.} We may proceed the other way round, starting with $g_1,\ldots,g_k$ instead of $G_1,\ldots,G_k$. Let $M_1,\ldots,M_k$ be finite non-empty subsets of $\mathbb{Z}^m$. We may fix $\cal F$ and $d_{ij}$ as before. Then we suppose that for $i=1,\ldots,k$ the following holds: $<p_j,q>\ge -d_{ij},j=1,\ldots,r$, for all $q\in M_i$. If $g_i\in\mathbb{C}[M_i]$, $i=1,\ldots,k$, i.e. $g_i=\sum_{q\in M_i} a_{iq} z^q$, we may pass from $g_1,\ldots,g_k$ to $G_1,\ldots,G_k$ and ${\rm g}_1,\ldots,{\rm g}_k$ as before.\\

 For $\sigma\in\cal F$, generated by $p_{j_1},\ldots,p_{j_s}$, let $M_i^\sigma:=\{q\in M_i\,|\,<p_{j_1},q>=-d_{ij_1},\ldots, <p_{j_s},q>=-d_{ij_s}\}$. Then put $g_i^\sigma:=\sum_{q\in M_i^\sigma} a_{iq} z^q$. So $g_i^\sigma\equiv 0$ if $M_i^\sigma=\emptyset$.\\

 We say that $\sigma$ is adapted to $M_1,\ldots,M_k$ and $(d_{ij})$ if the sets $M_1^\sigma,\ldots,M_k^\sigma$ are non-empty. If this holds for all $\sigma$ the fan $\cal F$ is called adapted to $M_1,\ldots,M_k$ and $(d_{ij})$. Obviously:\\

{\bf Lemma 1.4.1:} If $\sigma\in\cal F$ is not adapted to $M_1,\ldots,M_k$ and $(d_{ij})$ there is an $i$ such that $g_i^\sigma\equiv 0$.\\

By a suitable choice of $\cal F$ and $d_{ij}$ we can achieve that $\cal F$ is adapted to $M_1,\ldots,M_k$ and $(d_{ij})$:\\

 Let $\Delta$ be the convex hull of $M_1+\ldots+M_k$ in $\mathbb{R}^m$. If we assume that $\dim \Delta=m$, let ${\cal F}_0={\cal F}_0^+$ be the dual fan with respect to $\Delta$. If $\dim\Delta=m'$ is arbitrary we must be more careful: Fix $q\in\Delta$, let $M^+:=M\cap Span(\Delta-q)$, $N^+:=Hom(M^+,\mathbb{Z})$, then we have a fan ${\cal F}^+_0$ in $N^+_\mathbb{R}$ which is dual to $\Delta$. We have a canonical map $N_\mathbb{R}\longrightarrow N^+_\mathbb{R}$; taking the inverse images of  the cones of ${\cal F}^+_0$ we get a cone decomposition ${\cal F}_0$ of $N_\mathbb{R}$ which is only a fan if $m'=m$, we call it the dual cone decomposition. Let ${\cal F}_2$ be a fan which refines ${\cal F}_0$, it is complete; we have a toric morphism $X_{{\cal F}_2}\longrightarrow X_{{\cal F}^+_0}$. Let $p'_j, j=1,\ldots,r'$ be correspondingly defined (i.e. for ${\cal F}_2$ instead of ${\cal F}$), and let $d'_{ij}:= -\min\{<p'_j,q>\,|\,q\in M_i\}$, $i=1,\ldots,k$, $j=1,\ldots,r'$. Then we may apply the machinery above with $d'_{ij}$ instead of $d_{ij}$.\\
 
{\bf Lemma 1.4.2:}\\ 
a) The cones $\sigma\in{\cal F}$ which are adapted to $M_1,\ldots,M_k$ and $(d_{ij})$ form a subfan ${\cal F}_1$ of $\cal F$.\\
b) If $\cal F$ is adapted to $M_1,\ldots,M_k$ and $(d_{ij})$ the hypothesis of Lemma 1.2.2 holds, and $d_{ij}= -\min\{<p_j,q>\,|\,q\in M_i\}$.\\
c) With the notations above, ${\cal F}_2$ is adapted to $M_1,\ldots,M_k$ and $(d'_{ij})$.\\

 {\bf 1.5.} Now let us formulate non-degeneracy conditions: Let $M_1,\ldots,M_k,g_1,\ldots,g_k$ be chosen as in section 1.4.\\ 

 For $p\in \mathbb{R}^m$ let $M_i^p := \{q\in M_i\,|\,<p,q>=\min\{<p,q'>\,|\,q'\in M_i\}$, $g_i^p := \sum_{q\in M_i^p} a_{iq} z^q$. Then $g$ is called non-degenerate with respect to $M_1,\ldots,M_k$ if for every $p$ the mapping $(g_1^p,\ldots,g_k^p)$ $\colon (\mathbb{C}^*)^m\rightarrow \mathbb{C}^k$ has no critical point $z$ with $g_1^p(z)=\ldots=g_k^p(z)=0$. This condition is fulfilled if the coefficients $a_{iq}$ are chosen generically. If we call $g$ non-degenerate this is understood with respect to $supp\,g_1,\ldots, supp\,g_k$, where $supp\,g_i := \{q\,|\,a_{iq}\neq 0\}$ (under the hypothesis that $g_1\not\equiv 0,\ldots,g_k\not\equiv 0$). This is just the usual notion.\\

 On the other hand, let us call ${\rm g}:=({\rm g}_1,\ldots,{\rm g}_k)$ weakly non-degenerate (resp. non-degenerate) with respect to $M_1,\ldots,M_k$ if for every cone $\sigma\in{\cal F}$ the mapping $(g_1^\sigma,\ldots,g_k^\sigma)\colon (\mathbb{C}^*)^m\rightarrow \mathbb{C}^k$ has no critical point $z$ with $g_1^\sigma(z)=\ldots=g_k^\sigma(z)=0$ (resp. if moreover for every $p\notin |{\cal F}|$ the mapping $(g_1^p,\ldots,g_k^p)$ $\colon (\mathbb{C}^*)^m\rightarrow \mathbb{C}^k$ has no critical point $z$ with $g_1^p(z)=\ldots=g_k^p(z)=0$.) Then $Y$ is called a non-degenerate (resp. non-degenerate) complete intersection in $X$ with respect to $M_1,\ldots,M_k$. Of course the second condition can be dropped if $\cal F$ is complete.\\

 Let ${\cal F}_1$ be chosen as in Lemma 1.4.2a), $X_1:=X_{{\cal F}_1}$.\\

{\bf Lemma 1.5.1:} a) If $Y$ is a weakly non-degenerate complete intersection in $X$ with respect to $M_1,\ldots,M_k$ we must have that $Y$ is contained in $X_1$.\\
b) If $Y\subset X_1$, $Y$ is a weakly non-degenerate complete intersection in $X$ with respect to $M_1,\ldots$ $M_k$ if and only if 
$Y$ is a weakly non-degenerate complete intersection in $X_1$ with respect to $M_1,\ldots,M_k$.\\

 {\bf Proof:} a) If $\sigma$ is not adapted to $M_1,\ldots,M_k$ and $(d_{ij})$ we have $Y\cap O_\sigma=\emptyset$ because of Lemma 1.4.1.\\
 b) For $\sigma\in{\cal F}\setminus{\cal F}_1$ we have $Y\cap O_\sigma=\emptyset$, hence there is no critical point $z$ of $(g_1^\sigma,\ldots,g_k^\sigma):(\mathbb{C}^*)^m\to\mathbb{C}^k$ with $g_1^\sigma(z)=\ldots=g_k^\sigma(z)=0$.\\

We will now restrict to the case of an adapted fan.\\

{\bf Lemma 1.5.2:} Suppose that $\cal F$ is adapted to $M_1,\ldots,M_k$ and $(d_{ij})$. Then the following conditions are equivalent:\\
a) ${\rm g}$ is non-degenerate with respect to $M_1,\ldots,M_k$,\\
b) $g$ is non-degenerate with respect to $M_1,\ldots,M_k$.\\

 {\bf Proof:} If $\sigma\in{\cal F}, p\in \buildrel o\over\sigma$, we have $g^\sigma=g^p$.\\

 In this case $Y$ is a non-degenerate complete intersection with respect to $M_1,\ldots,M_k$ if the coefficients are chosen generically. Note that it is not necessary to specify the numbers $d_{ij}$ because we must have $d_{ij}=-\min\{<p_j,q>\,|\,q\in M_i\}$, by Lemma 1.4.2b).\\

{\bf Lemma 1.5.3:} Let $\rm g$ be weakly non-degenerate with respect to $M_1,\ldots,M_k$. Then the Jacobian matrix of $G$ along $Y'$ has rank $k$.\\

 {\bf Proof:} Use Lemma 1.3.4.\\

 Therefore we can apply Lemma 1.3.3. This means that every weakly non-degenerate complete intersection is quasi-smooth resp. smooth if $\cal F$ is simplicial resp. regular. Note that the stronger condition of being non-degenerate is useful in order to have a compactification with the same property:\\

{\bf Lemma 1.5.4:} Let $Y$ be a non-degenerate complete intersection with respect to $M_1,\ldots,$ $M_k$ and $\cal F$ adapted to $M_1,\ldots,M_k$ and $(d_{ij})$. Then there is a complete fan $\overline{\cal F}$ and a non-degenerate complete intersection $\overline{Y}$ in $\overline{X}:=X_{\overline{\cal F}}$ such that $\cal F$ is a subfan of $\overline{\cal F}$ and $Y=\overline{Y} \cap X$. If $\cal F$ is regular (resp. simplicial) $\overline{\cal F}$ can be chosen to be regular (resp. simplicial), too.\\

 {\bf Proof.} Let ${\cal F}^+_0$ and ${\cal F}_0$ be chosen as in the proof of Lemma 1.4.1. If we take the intersections of cones of ${\cal F}_0$ and $\cal F$ we just obtain the fan $\cal F$. So we take a corresponding suitable refinement $\overline{\cal F}$ of ${\cal F}_0$. (Note that we have a toric morphism $X_{\cal F}\longrightarrow X_{{\cal F}^+_0}$ which is compactified by $X_{\overline{\cal F}}\longrightarrow X_{{\cal F}^+_0}$.)\\

 So $Y$ is a Zariski open dense subset of some non-degenerate compact variety $\overline{Y}$ which can be chosen to be a (quasi-)smooth toric complete intersection if $Y$ is (quasi-)smooth.\\

 It is useful to introduce a weaker notion. Note that the cones of $\cal F$ which are contained in $\partial|{\cal F}|$ form a subfan $\partial{\cal F}$ of $\cal F$. In order to avoid complications we assume for the rest of section 1.5 that $\partial{\cal F}$ is adapted to $M_1,\ldots,M_k$ and $(d_{ij})$. Let us call $({\rm g}_1,\ldots,{\rm g}_k)$ non-degenerate at infinity with respect to $M_1,\ldots,M_k$ if the following holds:\par

 For every $p\notin |{\cal F}|$ the mapping $(g_1^p,\ldots,g_k^p)$ $\colon (\mathbb{C}^*)^m\rightarrow \mathbb{C}^k$ has no critical point $z$ with $g_1^p(z)=\ldots =g_k^p(z)=0$.\\

 It is easy to see that this condition is generically fulfilled. In this case let us call $Y$ non-degenerate at infinity with respect to $M_1,\ldots,M_k$. Of course the condition is automatically true if $\cal F$ is complete, that is why we speak of non-degeneracy at infinity. Furthermore $Y$ need no longer be a complete intersection. In the case $M_i=supp\,g_i,i=1,\ldots,k$, we simply say that $Y$ is non-degenerate at infinity. \\

 Then it is again possible to compactify in a suitable way. It is convenient to use in the case where $\cal F$ is complete the following notion: Let ${\cal F}_1$ be a subfan of $\cal F$ and $X_1:=X_{{\cal F}_1}$. Then $Y$ is called non-degenerate along $Y\setminus Y_1$ with respect to $M_1,\ldots,M_k$ if for all $\sigma\in{\cal F}\setminus{\cal F}_1$ the mapping $(g_1^\sigma,\ldots,g_k^\sigma)\colon (\mathbb{C}^*)^m\rightarrow \mathbb{C}^k$ has no critical point $z$ with $g_1^\sigma(z)=\ldots =g_k^\sigma(z)=0$.\\

{\bf Lemma 1.5.5:} Let $Y$ be non-degenerate at infinity with respect to $M_1,\ldots,M_k$. Then there is a complete fan $\overline{\cal F}$ and a subvariety $\overline{Y}$ of $\overline{X}:=X_{\overline{\cal F}}$ such that $\cal F$ is a subfan of $\overline{\cal F}$, $\overline{Y}\cap X=Y$ and $\overline{Y}$ is non-degenerate along $\overline{Y}\setminus Y$ with respect to $M_1,\ldots,M_k$. If $Y$ is (quasi-)smooth $\overline{Y}$ can be chosen to be (quasi-)smooth, too. Also, $Y$ is a complete intersection if and only if $\overline{Y}$ is a complete intersection.\\

 {\bf Proof.} Let us take the intersection of all halfspaces in $N_\mathbb{R}$ of the form $\{q\,|\,<p,q>\ge c\}$ with $p\in \overline{N_\mathbb{R}\setminus|{\cal F}|}$ which contain $M_1+\ldots+M_k$. Let ${\cal F}_0$ be the cone decomposition which is dual to this convex set, it is not necessarily complete. The cones contained in $\partial|{\cal F}|$
correspond to a subset ${\cal F}_1$ of ${\cal F}_0$. Note that $\partial{\cal F}$ (the subfan of $\cal F$ whose cones are contained in $\partial|{\cal F}|$ ) is a refinement of ${\cal F}_1$. So we can find a refinement of ${\cal F}_0$ which is a fan and contains $\partial{\cal F}$ as a subfan. If we add all cones of ${\cal F}\setminus\partial{\cal F}$ we obtain $\overline{\cal F}$. The rest is easy.\\ 

 Conversely, if $\overline{Y}$ is compact and non-degenerate along $\overline{Y}\setminus Y$ with respect to $M_1,\ldots,M_k$ we have that $Y$ is non-degenerate at infinity with respect to $M_1,\ldots,M_k$.\\

 {\bf 1.6.} We want to calculate certain invariants for smooth complete intersections which are non-degenerate at infinity with the property that they do not depend on the particular choice of the coefficients involved in the defining equations.\\

 Let us assume that $\cal F$ is regular and let us fix $d_{ij}$, $i=1,\ldots,k$, $j=1,\ldots,r$, and let $M_i\subset \{q\in\mathbb{Z}^m\,|\,<p_j,q>\ge -d_{ij},j=1,\ldots,r\}$. Put ${\bf M}:=M_1\times\cdots\times M_k$. Let $\partial{\cal F}$ be adapted to $M_1,\ldots,M_k$ and $d_{ij}$. Let us choose $\overline{\cal F}$ as in the proof of Lemma 1.5.5, $\overline{\cal F}$ regular. Then let $p_1,\ldots,p_{\bar{r}}$ be the generators of the corresponding edges, where $\bar{r}\ge r$, and $d_{ij}:=-\min\{<p_j,q>\,|\,q\in M_i\}$ for $j>r$. Let $X,\overline{X}, X',\overline{X}'$ be correspondingly defined. \\

 Put ${\cal X}':=X'\times \mathbb{C}^{\bf M}$, similarly ${\cal X},\overline{\cal X},\overline{\cal X}'$. Put ${\cal Y}':=\{(\zeta,z,(a_q))\in{\cal X}'\,|\,\sum_q a_{qi}\zeta^{\rho_i(q)}z^q=0, i=1,\ldots,k\}$, where $\rho_i(q):=(<p_1,q_i>+d_{i1},\ldots,<p_r,q_i>+d_{ir})$, similarly $\overline{\cal Y}'$. Furthermore, ${\cal Y}:=({\pi},id)({\cal Y}'), \overline{\cal Y}:=(\overline{\pi},id)(\overline{\cal Y}')$, where $\overline{\pi}:\overline{X}'\longrightarrow \overline{X}$ is defined in an obvious way. \\

 The canonical projection induces a mapping $p:{\cal X}\longrightarrow \mathbb{C}^{\bf M}$. For $(a_q)\in \mathbb{C}^{\bf M}$ we look at the following condition $(R)$:\\

$(R)$ The fibre of $p|{\cal Y}$ over $(a_q)$ is a smooth complete intersection which is non-degenerate at infinity with respect to $M_1,\ldots,M_k$.\\

 Let $S$ be the set of all points in $\mathbb{C}^{\bf M}$ where $(R)$ is not fulfilled. Then we have:\\

{\bf Theorem 1.6.1:} $S$ is a closed algebraic subset of $\mathbb{C}^{\bf M}$, and $p|{\cal Y}$ as well as $p|\overline{\cal Y}$ define topological fibre bundles over $\mathbb{C}^{\bf M}\setminus S$.\\

 {\bf Proof.} We can define a Whitney regular stratification of $\overline{X}$, taking $X$ and the orbits $O_\sigma, \sigma\in\overline{\cal F}\setminus{\cal F}$, as strata. Taking the product with $\mathbb{C}^{\bf M}$ we get a corresponding stratification of $\overline{\cal X}$. 
Now $\overline{\cal Y}$ is obtained by transversal intersection along $\overline{\cal X}\setminus{\cal X}$, in particular $\overline{\cal Y}$ is smooth along this set. 
So $\overline{\cal Y}\setminus {\cal Y}$ is endowed with a natural Whitney stratification. Now $S$ is the union of the critical values of the restriction of $p$ to ${\cal Y}$ and to the strata of $\overline{\cal Y}\setminus {\cal Y}$, hence a closed algebraic subset of $\mathbb{C}^{\bf M}$. Here a point of ${\cal Y}$ is called a critical point of $p|{\cal Y}$ if it is a critical point of $p|{\cal Y}_{reg}$ or does not belong to ${\cal Y}_{reg}$.
So we have a proper stratified submersion above the complement of $S$. Here we may apply Thom's first isotopy lemma in order to obtain topological fibre bundles.\\

 In particular, $\mathbb{C}^{\bf M}\setminus S$ is a Zariski-open subset of $\mathbb{C}^{\bf M}$, hence connected. So we get:\\

{\bf Corollary 1.6.2:} Given $M_1,\ldots,M_k$, ${\cal F}$ and $d_{ij}$, $\cal F$ adapted to $M_1,\ldots,M_k$ and $(d_{ij})$, all corresponding smooth complete intersections which are non-degenerate at infinity with respect to $M_1,\ldots,M_k$ are homeomorphic to each other.\\

 The varieties $Y$ for which we will give a method to compute the Hodge numbers are closed subvarieties of a compact toric variety which admit a partition into smooth locally closed subvarieties which are non-degenerate at infinity. The locally closed subvarieties are supposed to be the intersection of $Y$ by some $T$-invariant and irreducible locally closed subvariety of $X$; as we will see in section 1.7, such a subvariety of $X$ can be considered as a toric variety itself. An example of such a partition of $Y$ is given as follows:\\

{\bf Lemma 1.6.3:} Let $g_1,\ldots,g_k$ and $\cal F$ be as above, $\cal F$ complete. Suppose that for all $0<l\le k$, $1\le i_1<\ldots<i_l\le k$ the mapping $(g_{i_1},\ldots,g_{i_l})$ is non-degenerate with respect to $M_{i_1},\ldots,M_{i_l}$. Then the partition of $Y$ into the sets $O_\sigma\cap Y$, $\sigma\in {\cal F}$, is a partition into smooth locally closed subvarieties which are non-degenerate at infinity.\\

 If $\cal F$ is complete and the coefficients of $g_1,\ldots,g_k$ are chosen general enough $Y$ admits therefore 
a partition into smooth locally closed subvarieties which are non-degenerate at infinity.\\

 There are, however, interesting cases which are not covered by Lemma 1.6.3 directly but where it is possible to reduce to Lemma 1.6.3 by a homeomorphism, by Theorem 1.6.1:\\

{\bf Lemma 1.6.4}\\ 
a) Suppose that $Y$ is smooth and non-degenerate at infinity. By changing the coefficients of $g_1,\ldots,g_k$ we may obtain a homeomorphic variety $Y_1$ for which the intersections with all orbits are non-degenerate.\\
b) Suppose that $\cal F$ is complete and that $Y$ admits a decomposition into smooth complete intersections which are non-degenerate at infinity. Then the same conclusion as in a) holds.\\

 {\bf 2. Differential forms and Euler characteristics}\\

 {\bf 2.1.} Let $g_1,\ldots,g_k$ and $G_1,\ldots,G_k$ be chosen as in section 1.3. Let us suppose that ${\cal F}$ is complete. We assume that the condition (**) of the end of section 1.2 is fulfilled. Let $Y$ be correspondingly defined. We assume that $Y$ is a complete intersection, in fact it is sufficient to assume that $Y'$ is a complete intersection (see Lemma 1.3.2).\\

Before looking at the cohomology of differential forms let us consider the cohomology of $\pi_*{\cal O}_{X'}$. In particular we will compute $\chi(Y,{\cal O}_Y)$. Note that we can renounce here to the assumption that we have a simplicial fan!\\

 We can extend the action of $(\mathbb{C}^*)^r$ on $X'$ to an action of $(\mathbb{C}^*)^r\times(\mathbb{C}^*)^m$, where the action of $(\mathbb{C}^*)^m$ corresponds to the canonical action on $\mathbb{C}^m$. The corresponding  characters are given by $(s,q)\in \mathbb{Z}^r\times\mathbb{Z}^m$. This will lead to a finer graduation and a refined Poincar\'e series for ${\cal S}={\cal O}_{X'}$.\\

 Let $\sigma\in\cal F$ be generated by $p_1,\ldots,p_l$.\\

{\bf Lemma 2.1.1:} For $(s,q)\in \mathbb{Z}^r\times\mathbb{Z}^m$:
$\dim H^0(X_{\sigma'},{\cal O}_{X'})_{(s,q)}=1$ if $<p_j,q>+ s_j\ge 0$ for all $j\le l$,
$H^0(X_{\sigma'},{\cal O}_{X'})_{(s,q)}=0$ otherwise.\\

 {\bf Proof:} The only possible elements of $H^0(X_{\sigma'},{\cal O}_{X'})_{(s,q)}$ are of the form $c\zeta_1^{s_1+<p_1,q>}\cdot\ldots\cdot \zeta_r^{s_r+<p_r,q>}z^q, c\in\mathbb{C}$.\\

 For the higher cohomology groups we have of course, $X_{\sigma'}$ being affine:\\

{\bf Lemma 2.1.2:} $H^\lambda(X_{\sigma'},{\cal O}_{X'})=0$ for $\lambda>0$.\\

 Let us introduce the following formal Laurent series in the variables $x_1,\ldots,x_r,$ $t_1,\ldots,t_m$:
$$P_\sigma(x,t):=\sum \dim H^0(X_{\sigma'},{\cal O}_{X'})_{(s,q)} x^st^q$$
Let $\sigma_1,\ldots,\sigma_l$ be the maximal (i.e. $m$-dimensional) cones in $\cal F$. Then we have an open affine covering $X_{\sigma'_1},\ldots,$ $X_{\sigma'_l}$ of $X'$ which we can use in order to compute cohomology. \\

 In particular, putting $P(X',{\cal O}_{X'})(x,t):=\sum\chi(X',{\cal O}_{X'})_{s,q}x^st^q$,  we have\\
$$P(x,t):=P(X',{\cal O}_{X'})(x,t)=\sum_{1\le\nu\le l}\sum_{1\le j_1<\ldots<j_\nu\le l} (-1)^{\nu-1} P_{\sigma_{j_1}\cap\ldots\cap\sigma_{j_\nu}}(x,t)$$ 
Let $H(s,q)$ be the coefficient at $x^st^q$ in the series $P(x,t)$, i.e.
$$H(s,q):=\chi(X,(\pi_*{\cal O}_{X'})_{s,q})$$
 \\

 For $i=1,\ldots,l$ let $J_i\subset\{1,\ldots,r\}$ be defined as follows: $\sigma_i$ is generated by the $p_j$ with $j\in J_i$. Then $\sigma_{j_1}\cap\ldots\cap\sigma_{j_\nu}$ is generated by the $p_j$ with $j\in J_{j_1}\cap\ldots\cap J_{j_\nu}$. By Lemma 2.1.1 the coefficient of $P_{\sigma_{j_1}\cap\ldots\cap\sigma_{j_\nu}}(x,t)$ at $x^st^q$ is $1$ if $J_{j_1}\cap\ldots\cap J_{j_\nu}\subset I_{s,q}:=\{j\,|\,<p_j,q>\geq -s_j\}$ and $0$ otherwise.\\

 Now let $I$ be a subset of $\{1,\ldots,r\}$ and $\nu> 0$. Then let $m_{I,\nu}$ be the number of all subsets $K$ of $\{1,\ldots,l\}$ having $\nu$ elements , for which the intersection of all $J_\kappa$, $\kappa\in K$, is contained in $I$. Furthermore let $\chi_I:=\sum_{\nu>0} (-1)^{\nu-1} m_{I,\nu}$. Then:\\ 

{\bf Lemma 2.1.3:} \\
a) $H(s,q)=\chi(X,(\pi_*{\cal O}_{X'})_{s,q})=\chi_{I_{s,q}}$,\\
b) $P(x,t)=\sum_{s,q} \chi_{I_{s,q}} x^st^q$.\\

 Let $I$ be a subset of $\{1,\ldots,r\}$ and $s\in\mathbb{Z}^m$. Then let $n_{I,s}$ be the number of all $q\in\mathbb{Z}^m$ with $<p_j,q>\geq -s_j\Leftrightarrow j\in I$, $j=1,\ldots,r$, i.e. $n_{I,s}:=\{q\,|\,I=I_{s,q}\}$.\\

 Since $X$ is complete the vector spaces $H^\nu(X',{\cal O}_{X'})_s=H^\nu(X,(\pi_*{\cal O}_{X'})_s)$ are finite dimensional. Furthermore the dimension of $H^\nu(X,(\pi_*{\cal O}_{X'})_{s,q})$ depends only on $I_{s,q}$. Therefore $\chi_I\neq 0\Rightarrow n_{I,s}<\infty$. Altogether:\\

{\bf Lemma 2.1.4:}\\ 
a) $P(x):=P(x,1,\ldots,1)=\sum_s\sum_{I:\chi(I)\neq 0} \chi_I n_{I,s}x^s$,\\
b) $H(s):=\sum_q H(s,q)=\sum_s\sum_{I:\chi(I)\neq 0} \chi_In_{I,s}$.\\

{\bf Proposition 2.1.5:}\\ 
a) $\chi(Y,{\cal O}_Y)$ is equal to the coefficient at $x^0$ in the Laurent series $P(x)(1-x^{d_1})\cdots(1-x^{d_r})$,\\
b) $\chi(Y,{\cal O}_Y)=\sum_{\tau\ge 0,1\le l_1<\cdots<l_\tau\le k} (-1)^\tau H(-d_{l_1}-\ldots-d_{l_\tau})$.\\

{\bf Proof:} a) By induction on $j=0,\ldots,k$:\\
 $P(X',{\cal O}_{X'}/(G_1,\ldots,G_j))=P(x)(1-x^{d_1})\cdots(1-x^{d_j})$. Indeed, for $j=1,\ldots,k$ we have an exact sequence 
$$0\longrightarrow {\cal O}_{X'}/(G_1,\ldots,G_{j-1})\buildrel{\cdot G_j}\over\longrightarrow {\cal O}_{X'}/(G_1,\ldots,G_{j-1})\longrightarrow {\cal O}_{X'}/(G_1,\ldots,G_j)\longrightarrow 0$$
For $j=k$ we obtain the assertion.\\

 b) follows from a).\\

 In particular, this gives $\#Y$ if $\dim Y=0$ and the hypothesis of Lemma 1.3.3b) is fulfilled.\\

{\bf 2.2.} Now let us turn to differential forms. Therefore we want to work with manifolds or at least $V$-manifolds.\\

Therefore we suppose in this section from now on that the hypothesis of Lemma 1.3.3 is fulfilled.\\

In particular let $\cal F$ be a simplicial fan chosen as in \S 1, $\cal F$ complete. Let $Y\subset X_{\cal F}$ be accordingly defined. Then $Y$ is a compact $V$-manifold which is a complete intersection. \\

 We consider algebraic differential forms in the sense of Zariski, cf. [O]. First let us look at 1-forms. Let $\Omega_X^1$ be defined as follows: If $X$ is smooth it is defined as usual. If $X$ is quasi-smooth let $i\colon X_0\rightarrow X$ be the inclusion of the regular locus $X_0$, and let $\Omega_X^1 :=i_*\Omega^1_{X_0}$. Note that the holomorphic analogue would be the sheaf of weakly holomorphic 1-forms. The reason for this modification will be clear from Theorem 2.2.1.\\

 Similarly let $\Omega_{X'}^1$ be the sheaf of regular algebraic 1-forms on $X':=X_{\cal F'}$. The action of the torus $(\mathbb{C}^*)^r$ on $X'$ induces corresponding actions on the cohomology groups.\\

 Furthermore, the action of the torus leads to vector fields $D_1,\ldots,D_r$ on $X'$:
$$D_j := \zeta_j {\partial\over\partial\zeta_j} - p_{j1}z_1{\partial\over\partial z_1} -\ldots-p_{jm}z_m{\partial\over\partial z_m}$$
For $0\le\rho\le r$, we have the mapping $D^\rho:=(D_1,\ldots,D_\rho):\Omega^1_{Y'}\longrightarrow {\cal O}^\rho_{Y'}$. 
Let $\Omega^1_{Y',\rho}:=ker(D^\rho\colon \Omega^1_{Y'}\rightarrow {\cal O}_{Y'}^\rho)$.\\

{\bf Theorem 2.2.1:} $\Omega_{Y}^1 \simeq (\pi_*\Omega^1_{Y',r})_0$.\\

 {\bf Proof:} Let $Y'_0:=\pi^{-1}(Y_0)$ and $i'$ be the inclusion of $Y'_0$ in $Y'$, $\pi_0:=\pi|Y'_0:Y'_0\rightarrow Y_0$. Since $i\circ\pi_0=\pi\circ i'$ we have $i_*((\pi_0)_*\Omega^1_{Y'_0,r})=\pi_*(i'_*\Omega^1_{Y'_0,r})$. Now
$i'_*\Omega^1_{Y'_0,r}=i'_*(ker\,D^r:\Omega^1_{Y'_0}\rightarrow {\cal O}^r_{Y'_0})=\Omega^1_{Y',r}$ because $Y'$ is smooth, so $i_*((\pi_0)_*\Omega^1_{Y'_0,r})=\pi_*\Omega^1_{Y',r}$.
Taking care of the grading we get $i_*(((\pi_0)_*\Omega^1_{Y'_0,r})_0)=(\pi_*\Omega^1_{Y',r})_0$. Finally, 
$((\pi_0)_*\Omega^1_{Y'_0,r})_0=\Omega^1_{Y_0}$, so 
$\Omega^1_Y:=i_*\Omega^1_{Y_0}=i_*(((\pi_0)_*\Omega^1_{Y'_0,r})_0)=(\pi_*\Omega^1_{Y',r})_0$.
\\

{\bf Lemma 2.2.2:} For $\rho=1,\ldots,r$ we have an exact sequence
$$0\longrightarrow \Omega^1_{Y',\rho}\longrightarrow \Omega^1_{Y',\rho-1}\buildrel {D_\rho}\over\longrightarrow {\cal O}_{Y'}\longrightarrow 0$$\\

 {\bf Proof:} We need only show that the mapping induced by $D_\rho$ is surjective. Here it is sufficient to work with $X_{\sigma'}$ instead of $Y'$. Assume that $\sigma$ is generated by $p_{i_1},\ldots,p_{i_l}$. Let $h\in{\cal O}(X_{\sigma'})$ be given.\\

 If $\rho\notin \{i_1,\ldots,i_l\}$ we have that ${hd\zeta_\rho\over\zeta_\rho}\in \Omega^1_{X'}(X_{\sigma'})$ is mapped by $D_\rho$ onto $h$.\\

 So let us look at the case $\rho\in\{i_1,\ldots,i_l\}$. Note that $l\le m$ and there are $j_1,\ldots,j_l$ such that $det\,(p_{ij})_{i=i_1,\ldots,i_l,j=j_1,\ldots,j_l}\neq 0$. So we can find $\phi_i\in{\cal O}(X_{\sigma'}),i=1,\ldots,r, \psi_j\in{\cal O}(X_{\sigma'}),
j=1,\ldots,m$, such that $\phi_i=0$ for $i=i_1,\ldots,i_l$, $\phi_i-\sum_{j=1}^m p_{ij}\psi_j=0$ for $i=1,\ldots,r$, $i\neq\rho$ and $\phi_\rho-\sum_{j=1}^m p_{\rho j}\psi_j=h$. In fact, the $\psi_j,j\neq j_1,\ldots,j_l$ are arbitrary (e.g. $0$). Then, putting $\omega:=\sum_{i=1}^r{\phi_i\over\zeta_i}d\zeta_i+\sum_{j=1}^m{\psi_j\over z_j}dz_j$, we have $\omega\in \Omega^1_{X',\rho-1}(X_{\sigma'})$, $D_\rho(\omega)=h$.\\

 Note that we can associate to $\Omega^1_{X'}$ a Poincar\'e series: Let ${\cal S}$ be an arbitrary equivariant coherent sheaf on $X'$, for instance ${\cal S}=\Omega^1_{X'}$. Then $\pi_*{\cal S}$ is equipped with a grading: $\pi_*{\cal S}=\oplus_{s\in\mathbb{Z}^m} (\pi_*{\cal S})_s$. Now $\pi$ is affine, so $H^j(X,\pi_*{\cal S})=H^j(X',{\cal S})$. So the grading of $\pi_*{\cal S}$ induces a grading of $H^j(X',{\cal S})$. Since each $(\pi_*{\cal S})_s$ is a coherent ${\cal O}_X$-module the vector space $H^q(X',{\cal S})_s$ is finite dimensional and we can consider the Euler characteristic $\chi(X,{\cal S})_s:=\chi(X',(\pi_*{\cal S})_s)= \sum_q (-1)^q \dim H^q(X',{\cal S})_s$. So we can finally define the formal Laurent series
$$P(X',{\cal S})(x):=\sum_s \chi(X,{\cal S})_s x^s$$
In particular, we can look at ${\cal S}:=\Omega^1_{X'}$ . Then we will see:\\

{\bf Proposition 2.2.3:} $\chi(Y,\Omega^1_Y)$ is the coefficient at $x^0$ in the series $P(x)(1-x^{d_1})\ldots(1-x^{d_k})[x_1+\ldots+x_r+m-r-x^{d_1}-\ldots-x^{d_k}]$.\\

{\bf 2.3.} Now let us look at alternating differential forms: $\Omega^p_{Y'}:=\bigwedge^p \Omega^1_{Y'}$.\\

 Note that $\Omega^p_{X'}$ and $\Omega^p_{Y'}$ are equivariant coherent sheaves.\\

 Let $\Omega_{X'}:=\oplus_p \Omega^p_{X'}$ and $P(X,\pi_*\Omega_{X'})(x,t,y):= \sum_p P(X,\pi_*\Omega^p_{X'})(x,t)y^p$. The definition of $P(x,t)$ implies\\

{\bf Lemma 2.3.1:} $P(X',\Omega_{X'})(x,t,y):= \sum_p P(X',\Omega^p_{X'})(x,t)y^p=P(x,t)(1+yx_1)\ldots(1+yx_r)(1+y)^m$.\\

 {\bf Proof:} Note that $\Omega^1_{X'}$ is generated by $d\zeta_1,\ldots,d\zeta_r,{dz_1\over z_1},\ldots,{dz_m\over z_m}$.\\

 When we pass to a complete intersection we lose the grading with respect to $q$, so we put $t_1=\ldots=t_m=1$.\\

{\bf Lemma 2.3.2:} $P(Y',\Omega_{X'}/(G_1,\ldots,G_k))(x,y)= P(x)(1+yx_1)\ldots(1+yx_r)(1+y)^m(1-x^{d_1})\ldots(1-x^{d_k})$.\\

 {\bf Proof.} By induction on $j$ one shows for $j=1,\ldots,k$:\\
 $P(X',\Omega_{X'}/(G_1,\ldots,G_j))(x,y)= P(x)(1+yx_1)\ldots(1+yx_r)(1+y)^m(1-x^{d_1})\ldots(1-x^{d_j})$. This is because of the exact sequence
$$0\longrightarrow \Omega_{X'}/(G_1,\ldots,G_{j-1})\buildrel{\cdot G_j}\over\longrightarrow \Omega_{X'}/(G_1,\ldots,G_{j-1})\longrightarrow \Omega_{X'}/(G_1,\ldots,G_j)\longrightarrow 0$$
For $j=k$ we obtain the assertion.\\

{\bf Lemma 2.3.3:} $P(Y',\Omega_{Y'})(x,y)= P(x)(1+yx_1)\ldots(1+yx_r)(1+y)^m{1-x^{d_1}\over 1+yx^{d_1}}\ldots{1-x^{d_k}\over 1+yx^{d_k}}$.\\

 {\bf Proof.} For $j=1,\ldots,k$: $P(X',\Omega_{X'}/(G_1,\ldots,G_k,dG_1,\ldots,dG_j))(x,y)= P(x)(1+yx_1)$ $\ldots$ $(1+yx_r)(1+y)^m{(1-x^{d_1})\ldots(1-x^{d_k})\over (1+yx^{d_1})\ldots(1+yx^{d_j})}$. This is due to the exact sequence\\
 $0\longrightarrow \Omega^{p-1}_{X'}/(G_1,\ldots,G_k,dG_1,\ldots,dG_j)\buildrel{\wedge dG_j}\over\longrightarrow\Omega^p_{X'}/(G_1,\ldots,G_k,dG_1,\ldots,dG_{j-1})$\\
$\longrightarrow \Omega^p_{X'}/(G_1,\ldots,G_k,dG_1,\ldots,dG_j)\longrightarrow 0$\\
 For $j=k$ we obtain the assertion.\\

 For $\rho=0,\ldots,r$ let $\Omega^p_{Y',\rho} := \bigwedge^p\Omega^1_{Y',\rho}$, and let $\Omega^p_Y:=i_*\Omega^p_{Y_0}$.\\

{\bf Lemma 2.3.4:} For $\rho=0,\ldots,r$: $P(Y',\Omega_{Y',\rho})(x,y)= P(x)(1+yx_1)\ldots(1+yx_r)(1+y)^{m-\rho}{1-x^{d_1}\over 1+yx^{d_1}}\ldots{1-x^{d_k}\over 1+yx^{d_k}}$.\\

 {\bf Proof.} For $\rho >0$ we have an exact sequence 
$$0\rightarrow \Omega^p_{Y',\rho}\rightarrow \Omega^p_{Y',\rho-1}\buildrel D_\rho\over\rightarrow\Omega^{p-1}_{Y',\rho}\rightarrow 0$$
\\

 Finally:\\

{\bf Theorem 2.3.5:} $\chi(Y,\Omega^p_Y)$ is the coefficient at $x^0y^p$ in the series \\
${P(x)(1+yx_1)\ldots(1+yx_r)\over(1+y)^{r-m}}{1-x^{d_1}\over 1+yx^{d_1}}\ldots{1-x^{d_k}\over 1+yx^{d_k}}$.\\

 {\bf Proof.} We have $\Omega^p_Y=(\pi_*\Omega^p_{Y',r})_0$, as in Theorem 2.2.1.\\

In particular, we get Proposition 2.2.3.\\

 The use of the Poincar\'e series made a compact formula possible, for the actual computation the Hilbert function seems to be more useful:\\

{\bf Corollary 2.3.6:} $\chi(Y,\Omega^p_Y)=\sum_{\rho,\tau}(-1)^{p+\tau-\rho}\sum_{\rho+i_1+\cdots+i_k\le p} \sum_{1\le j_1<\cdots<j_\rho\le r}$ 
 $\sum_{1\le l_1<\cdots<l_\tau\le k} H(-e_{j_1}-\cdots-e_{j_\rho}-d_{l_1}-\cdots-d_{l_\tau}-i_1d_1-\cdots-i_kd_k)$.\\

 In the case $Y=X$ we can derive more exactly:\\

{\bf Theorem 2.3.7:} $\chi(X,\Omega^p_X)_q$ is the coefficient at $x^0y^pt^q$ in the series \\
$P(x,t)(1+yx_1)\ldots(1+yx_r)\over(1+y)^{r-m}$.\\

 {\bf 2.4.} Now let us look at symmetric instead of alternating differential forms. To indicate this we write $\hat{\Omega}$ instead of $\Omega$. So $\hat{\Omega}^p_{Y'} := S^p\Omega^1_{Y'}$ , where $S^p$ denotes the $p$-th symmetric tensor power.\\
\\

{\bf Lemma 2.4.1:} $P(X',\hat{\Omega}_{X'})(x,y,t)= {P(x,t)\over (1-yx_1)\ldots(1-yx_r)(1-y)^m}$.\\

{\bf Lemma 2.4.2:} $P(Y',\hat{\Omega}_{X'}/(G_1,\ldots,G_k))(x,y)= {P(x)(1-x^{d_1})\ldots(1-x^{d_k})\over(1-yx_1)\ldots(1-yx_r)(1-y)^m}$.\\

 {\bf Proof.} By induction on $j$ one shows for $j=1,\ldots,k$:\\
 $P(X',\hat{\Omega}_{X'}/(G_1,\ldots,G_j))(x,y)= {P(x)(1-x^{d_1})\ldots(1-x^{d_j})\over(1-yx_1)\ldots(1-yx_r)(1-y)^m}$, cf. Lemma 2.3.2. For $j=k$ we obtain the assertion.\\

{\bf Lemma 2.4.3:} $P(Y',\hat{\Omega}_{Y'})(x,y)=$  ${P(x)[(1-x^{d_1})(1-yx^{d_1})]\ldots[(1-x^{d_k})(1-yx^{d_k})]\over(1-yx_1)\ldots(1-yx_r)(1-y)^m}$.\\

 {\bf Proof.} For $j=1,\ldots,k$: $P(Y',\hat{\Omega}_{X'}/(G_1,\ldots,G_k,dG_1,\ldots,dG_j))(x,y)$\\
 $= {P(x)(1-x^{d_1})\ldots(1-x^{d_k})(1-yx^{d_1})\ldots(1-yx^{d_j})\over(1-yx_1)\ldots(1-yx_r)(1-y)^m}$. This is because of the exact sequence \\
 $0\longrightarrow \hat{\Omega}^{p-1}_{X'}/(G_1,\ldots,G_k,dG_1,\ldots,dG_{j-1})\buildrel{\otimes dG_j}\over\longrightarrow\hat{\Omega}^p_{X'}/(G_1,\ldots,G_k,dG_1,\ldots,dG_{j-1})$\\
$\longrightarrow \hat{\Omega}^p_{X'}/(G_1,\ldots,G_k,dG_1,\ldots,dG_j)\longrightarrow 0$\\
 For $j=k$ we obtain the assertion.\\

 For $\rho=0,\ldots,r$ let $\hat{\Omega}^p_{Y',\rho} := S^p\hat{\Omega}^1_{Y',\rho}$, $\hat{\Omega}^p_Y:=i_*\hat{\Omega}^p_{Y_0}$.\\

{\bf Lemma 2.4.4:} For $\rho\le r$: $P(Y',\hat{\Omega}_{Y',\rho})(x,y)= {P(x)(1-x^{d_1})(1-yx^{d_1})\ldots(1-x^{d_k})(1-yx^{d_k})\over(1-yx_1)\ldots(1-yx_r)(1-y)^{m-\rho}}$.\\

 {\bf Proof.} We have an exact sequence
$$0\rightarrow \hat{\Omega}_{Y',\rho}^p\rightarrow\hat{\Omega}_{Y',\rho-1}^p\buildrel D_\rho\over\rightarrow \hat{\Omega}_{Y',\rho-1}^{p-1}\rightarrow 0$$

 Finally:\\

{\bf Theorem 2.4.5:} $\chi(Y,\hat{\Omega}^p_Y)$ is the coefficient at $x^0y^p$ in the series\\
${P(x)(1-x^{d_1})(1-yx^{d_1})\ldots(1-x^{d_k})(1-yx^{d_k})(1-y)^{r-m}\over(1-yx_1)\ldots(1-yx_r)}$.
\\

 {\bf Proof.} We have $\hat{\Omega}_Y^p\simeq (\pi_*\hat{\Omega}_{Y',r}^p)_0$.\\

In particular, we get Proposition 2.2.3.\\

 In the case $Y=X$ we obtain more exactly:\\

{\bf Theorem 2.4.6:} $\chi(Y,\hat{\Omega}^p_Y)_q$ is the coefficient at $x^0y^pt^q$ in the series
${P(x,t)(1-y)^{r-m}\over(1-yx_1)\ldots(1-yx_r)}$.\\

 {\bf 2.5.} Furthermore let us look at differential forms without any symmetry condition. Instead of $\Omega$ we write $\check{\Omega}$ now, so 
$\check{\Omega}^p_{Y'}:=\otimes_{i=1}^p \Omega^1_{Y'}$,
$\check{\Omega}^p_Y:=i_*\check{\Omega}^p_{Y_0}$. First we have by induction:\\

{\bf Lemma 2.5.1:} 
$P(Y',\check{\Omega}^p_{Y',r})(x)= P(x)(1-x^{d_1})\ldots(1-x^{d_k})[x_1+\ldots+x_r+m-r-x^{d_1}-\ldots-x^{d_k}]^p$.\\

 For $p=0$ cf. Proposition 2.2.5. In order to prove the induction step we show:\\

{\bf Lemma 2.5.2:} 
$P(Y',\Omega^1_{X'}/(G_1,\ldots,G_k)\otimes_{{\cal O}_{Y'}}\check\Omega^p_{Y',r})(x)$   $=P(x)(1-x^{d_1})\ldots(1-x^{d_k})[x_1+\ldots+x_r+m-r-x^{d_1}-\ldots-x^{d_k}]^p(x_1+\ldots+x_r+m)$.\\

{\bf Lemma 2.5.3:} 
$P(Y',\Omega^1_{Y'}\otimes_{{\cal O}_{Y'}}\check{\Omega}^p_{Y',r})(x)=$   $P(x)(1-x^{d_1})\ldots(1-x^{d_k})[x_1+\ldots+x_r+m-r-x^{d_1}-\ldots-x^{d_k}]^p(x_1+\ldots+x_r+m-x^{d_1}-\ldots-x^{d_k})$.\\

 Finally we use the exact sequence $0\longrightarrow \Omega^1_{Y',r}\longrightarrow\Omega^1_{Y'}\buildrel D^r\over\longrightarrow {\cal O}_{Y'}^r\longrightarrow 0$.\\

 Lemma 2.5.1 implies for $\check{\Omega}_{Y',r}:=\oplus_p\check{\Omega}^p_{Y',r}$ :\\

{\bf Lemma 2.5.4:} $P(Y',\check{\Omega}_{Y',r})(x,y):=\sum_p P(Y',\check{\Omega}^p_{Y',r})(x)y^p=$\\
 ${P(x)(1-x^{d_1})\ldots(1-x^{d_k})\over 1-y(x_1+\ldots+x_r+m-r-x^{d_1}-\ldots-x^{d_k})}$.\\

{\bf Theorem 2.5.5:} $\chi(Y,\check{\Omega}^p_Y)$ is the coefficient at $x^0y^p$ in\\  
${P(x)(1-x^{d_1})\ldots(1-x^{d_k})\over 1-y(x_1+\ldots+x_r+m-r-x^{d_1}-\ldots-x^{d_k})}$.\\

 {\bf Proof.} We have $\check{\Omega}^p_Y=(\pi_*\check{\Omega}^p_{Y',r})_0$.\\

In particular, we get Proposition 2.2.3 again.\\

 In the case $Y=X$ we have an additional grading and obtain in the same way as before:\\

{\bf Theorem 2.5.6:} $\chi(X,\check{\Omega}^p_X)_q$ is the coefficient at $x^0y^pt^q$ in  ${P(x,t)\over 1-y(x_1+\ldots+x_r+m-r)}$.
\\

 {\bf 2.6.} As an example let us take the case of complete intersections in weighted projective spaces. Let $w_1,\ldots,w_m$ be positive integers which are relatively prime. Then we have a corresponding grading for $R:=Spec\,\mathbb{C}[z_1,\ldots,z_m]$: $R=\oplus_{t\ge 0} R_t$, where $R_t$ is spanned by the monomials $z_1^{j_1}\cdots z_m^{j_m}$ with $w_1j_1+\ldots+w_mj_m=t$. Then $\mathbb{P}_{(w_1,\ldots,w_m)}=Proj\,R$ is the corresponding weighted projective space. The transcendental topology gives $\mathbb{P}_{(w_1,\ldots,w_m)}\simeq (\mathbb{C}^m\setminus\{0\})/\sim$, where $z\sim z'\Leftrightarrow$ there is a $c\in\mathbb{C}^*$ such that $z'_j=c^{w_j}z_j$, $j=1,\ldots,m$.\\

 In fact, weighted projective spaces are toric varieties. It is easy to proceed as follows: Let $w_0:=1$. Then $\mathbb{P}_{(w_0,\ldots,w_m)}$ can be identified with $X=X_{\cal F}$ where ${\cal F}$ consists of all cones in $\mathbb{R}^m$ which are generated by at most $m$ of the vectors $p_0,\ldots,p_m$, where $p_0=(-w_1,\ldots,-w_m)$ and $p_j=e_j$, $j=1,\ldots,m$. The fan is complete and simplicial but not necessarily regular, so $X$ is compact and quasi-smooth. Furthermore, let $\sigma^0$ be the cone spanned by $p_0$, then $X\cap \overline{O_{\sigma^0}}\simeq \mathbb{P}_{(w_1,\ldots,w_m)}$.\\

 Let $P(x)$ be defined as in section 2.1. We want to calculate $P(x)$. Let $\sigma_j$ be the cone generated by $p_0,\ldots,p_{j-1},p_{j+1},\ldots,p_m$, so $J_j=\{0,\ldots,j-1,j+1,\ldots,m\}$. If $I$ is a subset of $\{0,\ldots,m\}$ having $\mu$ elements we have $J_{j_1}\cap\ldots\cap J_{j_\nu}\subset I$ if and only if $\{0,\ldots,m\}\setminus \{j_1,\ldots,j_\nu\}$ is contained in $I$. Therefore $m_{I,\nu}={\mu\choose m+1-\nu}$. Putting $\lambda:=m+1-\nu$ we get $\chi_I=\sum_{\lambda=0}^m(-1)^{m-\lambda}{\mu\choose\lambda}$. So $\chi_I=1$ for $\mu=m+1$, $\chi_I=(-1)^m$ for $\mu=0$ and $\chi_I=0$ for $0<\mu\le m$.\\

 Now the number of all $q$ such that $<p_j,q>\ge -s_j, j=0,\ldots,m$ is by definition equal to $n_{\{0,\ldots,m\},s}$. We have the following formula:\\

{\bf Lemma 2.6.1:} The number of all $q$ such that $<p_j,q>\ge -s_j, j=0,\ldots,m$ is\\
$n_{\{0,\ldots,m\},s}=res_0 {x^{-s_0-1}\over 1-x}{x^{-w_1s_1}\over 1-x^{w_1}}\cdots{x^{-w_ms_m}\over 1-x^{w_m}}$.\\

{\bf Proof:} We develop the function which appears on the right hand side in a Laurent series: $\sum x^{-1+q_0+q_1w_1+\ldots+q_mw_m}$, where the sum extends over all $(q_0,\ldots,q_m)\in\mathbb{Z}^{m+1}$ such that $q_j\ge -s_j, j=0,\ldots,m$. So $res_0=\#\{(q_0,\ldots,q_m)\,|\,q_j\ge-s_j, j=0,\ldots,m, q_0+q_1w_1+\ldots+q_mw_m=0\}=\#\{(q_1,\ldots,q_m)\,|\,q_j\ge-s_j, j=1,\ldots,m, q_1w_1+\ldots+q_mw_m\le s_0\}$, which is the number on the left hand side.\\

On the other hand, the number of all $q$ such that $<p_j,q>< -s_j, j=0,\ldots,m$ is equal to the number of all $q'$ such that $<p_j,q'>\ge s_j+1,j=0,\ldots,m$, i.e.
$$n_{\emptyset,s}= res_0 {x^{s_0}\over 1-x}{x^{w_1(s_1+1)}\over 1-x^{w_1}}\cdots{x^{w_m(s_m+1)}\over 1-x^{w_m}} =(-1)^m res_\infty {x^{-s_0-1}\over 1-x}{x^{-w_1s_1}\over 1-x^{w_1}}\cdots{x^{-w_ms_m}\over 1-x^{w_m}}$$ 
as one sees using the substitution $\xi=x^{-1}$ and Lemma 2.6.1. So we get, with $res_{0,\infty}:=res_0+res_\infty$:\\

{\bf Lemma 2.6.2:} $P(x_0,\ldots,x_m)=\sum H(s_0,\ldots,s_m)x_0^{s_0}\cdots x_m^{s_m}$ with 
$H(s_0,\ldots,s_m)=res_{0,\infty} {x^{-s_0-1}\over 1-x}{x^{-w_1s_1}\over 1-x^{w_1}}\cdots{x^{-w_ms_m}\over 1-x^{w_m}}$.\\

The following result will be useful:\\

{\bf Lemma 2.6.3:} If $Q(x_0,\ldots,x_m)$ is a Laurent series, the coefficient in $P\cdot Q$ at $x_0^{s_0}\cdot\ldots\cdot x_m^{s_m}$ is:\\
$res_{0,\infty} {x^{-s_0-1}\over 1-x}{x^{-w_1s_1}\over 1-x^{w_1}}\cdots{x^{-w_ms_m}\over 1-x^{w_m}}Q(x,x^{w_1},\ldots,x^{w_m})$.\\

{\bf Proof:} Put $Q=\sum_sQ_sx^s$. Then the coefficient in $P\cdot Q$ at $x_0^{s_0}\cdot\ldots\cdot x_m^{s_m}$ is:\\
$\sum_t res_{0,\infty} {x^{-t_0-1}\over 1-x}{x^{-w_1t_1}\over 1-x^{w_1}}\cdots{x^{-w_mt_m}\over 1-x^{w_m}}Q_{s-t}$\\
$=res_{0,\infty} {x^{-s_0-1}\over 1-x}{x^{-w_1s_1}\over 1-x^{w_1}}\cdots{x^{-w_ms_m}\over 1-x^{w_m}}\sum_tQ_{s-t}x^{(s_0-t_0)+w_1(s_1-t_1)+\ldots+w_m(s_m-t_m)}$.\\
Furthermore, $\sum_tQ_{s-t}x^{(s_0-t_0)+w_1(s_1-t_1)+\ldots+w_m(s_m-t_m)}=Q(x,x^{w_1},\ldots,x^{w_m})$.\\

 Now let $\tilde{g}_1,\ldots,\tilde{g}_k\in\mathbb{C}[z_0,\ldots,z_m]$ be weighted homogeneous polynomials of degree $d_1,\ldots,$ $d_k$ with respect to the weights $w_0,\ldots,w_m$, i.e. $\tilde{g}_i$ is a linear combination of monomials $z_0^{j_0}\ldots z_m^{j_m}$ with $w_0j_0+\ldots+w_mj_m=d_i$. If we put $g_i(z_1,\ldots,z_m):=\tilde{g}_i(1,z_1,\ldots,z_m)$, $d_{0i}:=d_i$ and $d_{li}:=0, l=1,\ldots,m$, we can apply section 1. Note that $Y$ is a quasi-smooth complete intersection as soon as $\tilde{g}_1,\ldots,\tilde{g}_k$ define a complete intersection with an isolated singularity. By the previous results we get:\\

{\bf Theorem 2.6.4:}\\
a) $\chi(Y,\Omega^p_Y)= res_{x=0,\infty} res_{y=0} {x^{-1}y^{-p-1}\over 1+y}{1+yx^{w_0}\over 1-x^{w_0}}\cdots{1+yx^{w_m}\over 1-x^{w_m}}{1-x^{d_1}\over 1+yx^{d_1}}\cdots {1-x^{d_k}\over 1+yx^{d_k}}$\\
b) $\chi(Y,\hat{\Omega}^p_Y)=res_{x=0,\infty} res_{y=0}{x^{-1}y^{-p-1}(1-x^{d_1})(1-yx^{d_1})\cdots (1-x^{d_k})(1-yx^{d_k})(1-y)\over (1-x^{w_0})(1-yx^{w_0})\cdots(1-x^{w_m})(1-yx^{w_m})}$\\
c) $\chi(Y,\check{\Omega}^p_Y)=res_{x=0,\infty} res_{y=0}{x^{-1}y^{-p-1}(1-x^{d_1})\cdots(1-x^{d_k})\over (1-x^{w_0})\cdots(1-x^{w_m})(1-y(x^{w_0}+\cdots+x^{w_m}-x^{d_1}-\cdots-x^{d_k}-1))}$\\

 We can also treat complete intersections in quasi-projective spaces, passing to $Y\cap\overline{O_{\sigma^0}}$: look at the complete intersection defined by $z_0,g_1,\ldots,g_k$. In this way it is easy to see that Theorem 2.6.4 holds without the assumption that $w_0=1$ provided that $w_0,\ldots,w_m$ are relatively prime. In particular part a) proves then a formula which was announced in [H1].\\

 {\bf 3. Computation of Hodge numbers}\\

 {\bf 3.1.} Recall that the cohomology of an algebraic variety $Y$ is equipped with a canonical mixed Hodge structure, according to Deligne [D]: On $H^j(Y;\mathbb{C})$ we have an increasing filtration $W$ and a decreasing filtration $F$ ; let $h^{jpq}(Y):=\dim_\mathbb{C} Gr^W_{p+q} Gr_F^p H^j(Y;\mathbb{C})$. If $Y$ is a compact $V$-manifold we have $h^{jpq}(Y)=0$ if $j\neq p+q$, so let $h^{pq}(Y):=h^{p+q,p,q}(Y)$ in this case. If $Y$ is even smooth and projective this coincides with the classical concept of Hodge numbers.\\

 The Euler-Hodge characteristics are the numbers $e^{pq}(Y):=\sum_j (-1)^j h^{jpq}(Y)$. In the case of a compact $V$-manifold it is sufficient to calculate these numbers, because $e^{pq}=(-1)^{p+q}h^{pq}$. If we work with cohomology with compact supports we write $e^{pq}_c(Y)$. Note that $e^{pq}_c$ is called $e^{pq}$ in [D-K].\\

 The aim of section 3 is to compute the Hodge numbers $h^{pq}$ for compact quasi-smooth varieties $Y$ which admit a decomposition into smooth complete intersections which are non-degenerate at infinity.\\

 As we will see, the essential step is to calculate the numbers $e^{pq}$ for non-degenerate complete intersections in tori. But first we argue as follows:\\

 If $Y$ is decomposed into $Y^i$, $i=1,\ldots,l$, we have $e^{pq}(Y)=e^{pq}_c(Y)=\sum_{i=1}^l e^{pq}_c(Y^i)$. So it is sufficient to calculate the numbers $e^{pq}_c$ for smooth toric complete intersections which are non-degenerate at infinity.\\

{\bf Theorem 3.1.1:} The numbers $e^{pq}_c$ of a smooth complete intersection $Y$ which is non-degenerate at infinity do not depend on the particular choice of the coefficients.\\

 {\bf Proof.} We use the results of sections 1.5 and 1.6 and induction on $\dim Y$. By Lemma 1.5.5 we have a smooth compactification $\overline{Y}$ of $Y$ which is non-degenerate along $\overline{Y}\setminus Y$. Now $h^{pq}(\overline{Y})=\dim H^q(\overline{Y},\Omega^p_{\overline{Y}})$. So $h^{pq}(\overline{Y})$ depends on the coefficients of the equations involved in an upper semicontinuous way. Since $\sum_{p+q=r} h^{pq}(\overline{Y})=b^r(\overline{Y})=$ $r$-th Betti number and the Betti numbers are constant by Theorem 1.6.1 we get that the numbers $h^{pq}(\overline{Y})$ do not depend on the particular choice of the coefficients. An analogous statement holds for the subsets of $\overline{Y}\setminus Y$ which are of the form $\overline{Y}\cap \overline{O_\sigma}$, by induction. By the additivity of the numbers $e^{pq}_c$ we get the assertion.\\

 So it is sufficient to calculate the numbers $e^{pq}(Y)$, $Y$ being defined as in section 1, provided that the coefficients involved are chosen generically. This means that we fix $M_1,\ldots,M_k$ and choose the coefficients $a_{iq}$ generically. It is now no longer important that $Y$ should be a smooth complete intersection which is non-degenerate at infinity.\\

 Now $Y$ decomposes into the sets $Y\cap O_\sigma$, so it is sufficient to compute the numbers $e^{pq}_c(Y\cap O_\sigma)$.\\

 Since the coefficients are chosen generically $Y\cap O_\sigma$ is a smooth complete intersection in $O_\sigma$ - but not necessarily of codimension $k$ ! This is because some of the sets $M_1^\sigma,\ldots,M_k^\sigma$ may be empty. Anyhow we are left with the question how to compute the numbers $e^{pq}_c(Y^*)$ where $Y^*=\{z\in(\mathbb{C}^*)^m\,|\,g_1(z)=\ldots=g_k(z)=0\}$ (perhaps with a different $k$) and the coefficients of $g_1\not\equiv 0,\ldots,g_k\not\equiv 0$ are chosen generically. For the sake of completeness we note\\

{\bf Corollary 3.1.2:} The numbers $e^{pq}(Y^*)$ depend only on the supports of $g_1,\ldots,g_k$ if $g=(g_1,\ldots,g_k)$ is non-degenerate.\\

 So it is sufficient to compute the numbers $e^{pq}(Y^*)$.\\

 {\bf 3.2.} Now let us compute the numbers $e^{pq}_c$ for $Y^*=\{z\in(\mathbb{C}^*)^m\,|\,g_1(z)=\ldots=g_k(z)=0\}$ where $g=(g_1,\ldots,g_k)$ is non-degenerate. Because of section 3.1, we may suppose that for all $1\le i_1<\ldots<i_s\le k$ the mapping $(g_{i_1},\ldots,g_{i_s})$ is non-degenerate, too. In fact, the proof in [D-K] needs this assumption too, without being mentioned explicitly.\\

 By induction on $m$ and - for fixed $m$ - on $k$ let us calculate the numbers $e^{pq}(Y^*)$.\\

 Induction step: We distinguish the cases $\dim\Delta<m$ and $\dim\Delta=m$, where $\Delta$ is the convex hull of $supp\,g_1+\ldots+supp\,g_k$.\\

 If $\dim\Delta<m$, we may reduce to the case of a torus of smaller dimension, using K\"unneth formula.\\

 So suppose that $\dim\Delta=m$. Let ${\cal F}_0$ be the dual fan to $\Delta$, $\tilde{Y}=\{z\in(\mathbb{C}^*)^m\,|\,g_1(z)\cdot\ldots\cdot g_k(z)=0\}$ Then we have the following Lefschetz theorem:\\

{\bf Proposition 3.2.1:} The pair $((\mathbb{C}^*)^m,\tilde{Y})$ is $(m-1)$-connected.\\

 {\bf Proof.} There is an ample sheaf $\cal S$ on $X_{{\cal F}_0}$ associated with $\Delta$, see [O] Theorem 2.22, and  $g_1\cdot\ldots\cdot g_k$ can be interpreted as a section of $\cal S$, so $X_{{\cal F}_0}\setminus\{g_1\cdot\ldots\cdot g_k=0\}$ is affine. Therefore there is a fundamental system of neighbourhoods $U$ of $\{g_1\cdot\ldots\cdot g_k=0\}$ in $X_{{\cal F}_0}$ such that the pair $((\mathbb{C}^*)^m,(\mathbb{C}^*)^m\cap U)$ is
$(m-1)$-connected. Because of the nondegeneracy assumption $\tilde{Y}$ is a deformation retract of $(\mathbb{C}^*)^m\cap U$. So we get the assertion.\\

 There are other versions of Lefschetz theorems for the torus, see [Ok].\\

 So we get:\\

{\bf Lemma 3.2.2:} $e^{pq}(\tilde{Y}) = e^{pq}((\mathbb{C}^*)^m)$ f\"ur $p+q<n:=m-k$.\\

 {\bf Proof.} Assume that $p+q<n$. Since $g_1\cdot\ldots\cdot g_k=0$ defines a divisor with normal crossings in $(\mathbb{C}^*)^m$, $Gr^W_lH^j(\tilde{Y};\mathbb{Q})=0$ for $l\le j-k$, so $h^{jpq}(\tilde{Y})=h^{jpq}((\mathbb{C}^*)^m)=0$ for $j\geq m-1$. For $j<m-1$ we may apply Proposition 3.2.1 and obtain that $h^{jpq}(\tilde{Y})=h^{jpq}((\mathbb{C}^*)^m)$.\\ 

 Now we can compute the cohomology of $\tilde{Y}$ by a spectral sequence which involves the cohomology of the spaces $(\mathbb{C}^*)^m\cap\{g_{i_1}=\ldots= g_{i_s}=0\}$ with $1\leq i_1<\ldots<i_s \leq k$. In particular we obtain\\

{\bf Lemma 3.2.3:} $e^{pq}(\tilde{Y})=\sum_{s>0, 1\le i_1<\ldots<i_s\le k} (-1)^{s-1}e^{pq}((\mathbb{C}^*)^m\cap\{g_{i_1}=\ldots= g_{i_s}=0\})$.\\

 By induction the numbers $e^{pq}((\mathbb{C}^*)^m\cap\{g_{i_1}=\ldots=g_{i_s}=0\})$ with $s<k$ are known. Using Lemma 3.2.2 we obtain the numbers $e^{pq}(Y^*)$ with $p+q<n$.\\

 By duality we get $e^{pq}_c(Y^*)$ for $p+q>n$.  \\

 Let $\cal F$ a simplicial subdivision of ${\cal F}_0$. Let $Y$ be correspondingly defined; $Y$ is a $V$-manifold. By induction hypothesis $e^{pq}_c(Y\setminus Y^*)$ is known, so $e^{pq}(Y)$ for $p+q>n$. By Poincar\'e duality $e^{pq}(Y)$ is known for $p+q<n$, too. What is missing is the computation of $e^{pq}(Y)$ for $p+q=n$, which is accomplished with the help of the numbers $e^p(Y)=\sum_q e^{pq}(Y)=\chi(X,\Omega^p_Y)$ calculated in Theorem 2.3.5. So in total we know the numbers $e^{pq}(Y)$.\\

 As said before, the numbers $e^{pq}_c(Y\setminus Y^*)$ are known, so all $e^{pq}_c(Y^*)$, too. \\

 {\bf Example 3.2.4:} Let $\cal F$ be the fan which consists of all cones which are spanned by at most four of the vectors $e_1,e_2,e_3,e_4,-e_1-\ldots-e_4$ in $\mathbb{R}^4$. Then $X\simeq \mathbb{P}_4$. Let $g_1(z):=z_1^2+z_2^2+z_3^2+z_4^2-1$, $g_2(z):=z_1^3-z_2^3+z_3^3-z_4^3$.\\
 Then $g_1,g_2$ define a smooth complete intersection in $X$, and $g_1$, $g_2$ resp. $g_1g_2$ define hypersurfaces $Y_1, Y_2$ resp. $\tilde{Y}$. Note that not all of the hypersurfaces  are smooth.\\
 First we calculate the numbers $e^{pq}$ for the intersection of the hypersurfaces by $T$, $p+q<2=\dim\,Y$.\\
 By Lemma 3.2.2, we have $e^{pq}(Y_1\cap T)=e^{pq}(\tilde{Y}\cap T)=e^{pq}(T)$ for $p+q<2$.\\
 The support of $g_2$, however, spans a hypersurface, so $Y_2\cap T$ is a total space of some $\mathbb{C}^*$ bundle, namely with base space $B:=\{z\in(\mathbb{C}^*)^3\,|\,z_1^3-z_2^3+z_3^3=1\}$. Again, we have $e^{pq}(B)=e^{pq}((\mathbb{C}^*)^3)$ for $p+q<2$. By K\"unneth we get $e^{pq}(Y_2\cap T)$, $p+q<2$.\\
 Now we obtain the numbers $e^{pq}(Y\cap T)=e^{pq}(Y_1\cap T)+e^{pq}(Y_2\cap T)-e^{pq}(\tilde{Y}\cap T)$ for $p+q<2$.\\
 By duality we get the numbers $e^{pq}_c(Y\cap T)$, $p+q>2$.\\
 Now $e^{pq}(Y)=e^{pq}_c(Y\cap T)+ e^{pq}(Y\cap(X\setminus T))$. Since $\dim\,Y\cap(X\setminus T)\le 1$ we have that the numbers $e^{pq}$ of this space with $p+q>2$ vanish, so the numbers $e^{pq}(Y)=e^{pq}_c(Y\cap T)$ for $p+q>2$ are known.\\
 By duality we know $e^{pq}(Y)$ for $p+q<2$, too. So it is sufficient to calculate $e^p(Y)$, where we can use Theorem 2.6.4.\\

{\bf Example 3.2.5.}  Let $N=\mathbb{Z}^3$ and $\cal F$ be the fan in $\mathbb{R}^3$ whose cones are generated by at most three of the following vectors: $p_0=(-w_1,-w_2,-w_3)$, $p_j=e_j$, $j=1,2,3$ with $w_1=4, w_2=2,w_3=3$. Then $X:=X_{\cal F}$ is the weighted projective space $\mathbb{P}_w$, with $w=(w_0,\ldots,w_3)$ and $w_0=1$. It is quasi-smooth but not smooth because $\cal F$ is simplicial but $p_0,p_3$ cannot be extended to a basis of $\mathbb{Z}^3$. Note that $\mathbb{C}^3$ can be regarded as an open subset of $T$, it corresponds to the cone generated by $p_1,p_2,p_3$. See section 2.6.\\
Now let $g_0:=z_1^3+z_2^6+z_3^4$ and $g:=g_0-2z_3^2+z_2+1$. Note that $g_0$ is a convenient weighted homogeneous polynomial (with respect to the weights above) with an isolated singularity. \\
By Lemma 1.3.3 it is easy to see that $g$ defines a quasi-smooth hypersurface $Y$ in $X$. Note that $Y$ intersects the singular locus of $X$.\\
On the other hand, $g$ is not non-degenerate: For $p:=(1,1,0)$, we have $g^p=(z_3^2-1)^2$. So we cannot work with $Y\cap T$ directly. In fact, $Y$ is tangent to the orbit $\{0\}\times\mathbb{C}^*$ at the two points $(0,0,\pm 1)$.\\
But $Y\cap\mathbb{C}^3$ is smooth and non-degenerate at infinity. The Euler-Hodge characteristics of this space do not change if we pass from $g$ to $g_1:=g_0-1$. So we can apply the procedure above. \\
Let $Y$ now be defined using $g_1$. Then $e^{pq}(Y\cap T)=e^{pq}(T)$, $p+q<2$, so these numbers are known. So we get $e^{pq}_c(Y\cap T)$, $p+q>2$. Similarly as in the preceding example we have $e^{pq}(Y)=e^{pq}_c(Y\cap T)$ for these $p,q$. By duality we know $e^{pq}(Y)$ for $p+q<2$, too. Since we can calculate $e^p(Y)$ we obtain all numbers $h^{pq}(Y)$.\\

 {\bf Example 3.2.6.} Let $N=\mathbb{Z}^4$ and $\cal F$ be the fan in $\mathbb{R}^4$ whose cones are generated by at most four of the following vectors: $p_0=(-w_1,-w_2,-w_3,-w_4)$, $p_j=e_j$, $j=1,\ldots,4$ with $w_1=w_2=w_3=2, w_4=1$. Then $X:=X_{\cal F}$ is the weighted projective space $\mathbb{P}_w$, with $w=(w_0,\ldots,w_4)$ and $w_0=1$. Let $g_0:=z_1^2+z_2^2+z_3^2+z_1z_4^2$, $g:=g_0-1$. Note that $g_0$ belongs to Orlik's list of non-convenient weighted homogeneous polynomials (with respect to the weights above) with an isolated singularity.\\
By Lemma 1.3.3 it is again easy to see that $g$ defines a quasi-smooth hypersurface $Y$ in $X$. Note that $Y$ intersects the singular locus $\Sigma$ of $X$. In weighted homogeneous coordinates, the latter is given by $z_0=z_4=0$.\\
Here the intersection of $Y$ with an orbit of $X$ may be an entire orbit instead of a hypersurface in the orbit: $z_0=z_1=z_2=z_4=0$. But $Y\setminus\Sigma$ is non-degenerate at infinity. Therefore we can proceed as in the preceding example.\\
Note, however, that we can no longer expect that $e^{pq}(Y)=e^{pq}_c(Y\cap T)$, $p+q>\dim\,Y=3$, because we have to take $e^{22}_c(Y\cap O_\sigma)$ into account if $\dim\,Y\cap O_\sigma= 2$.\\

{\bf 3.3.} As we have seen in the examples the method described above can lead to complicated calculations. In special cases, however, shortcuts can reduce the labour consideraby.\\

In fact, the examples above belong to a class of examples which can be treated much more easily: complete intersections in weighted projective cases. The results have already been announced in  [H1]. They generalize the calculation made by Hirzebruch [Hi2] Theorem 22.1.1 for complete intersections in projective space - in fact he assumed that the participating hypersurfaces are smooth, too, which is not true in Example 3.2.4. Numerically, however, this plays no role, as one sees by some deformation argument.\\

Let $X:=\mathbb{P}_w$, $w=(w_0,\ldots,w_m)$. Let $Y$ be a complete intersection in $X$ of dimension $n$. Then we have a Lefschetz theorem for rational cohomology:\\

{\bf Lemma 3.3.1:} The mapping $H^j(X;\mathbb{Q})\to H^j(Y;\mathbb{Q})$ is an isomorphism for $j<n$ and injective for $j=n$.\\

{\bf Proof:} By duality we have $H^j(X,Y;\mathbb{Q})\simeq H_{2m-j}(X\setminus Y;\mathbb{Q})$. So we must show that $H^j(X\setminus Y;\mathbb{Q})=0$ for $j\ge m+k$. But we have $Y=Y_1\cap\ldots\cap Y_k$, where $Y_i$ is defined by $g_i=0$, and $X\setminus Y_{i_1}\cup\ldots\cup Y_{i_l}$ is affine, so $H^j(X\setminus Y_{i_1}\cup\ldots\cup Y_{i_l};\mathbb{Q})=0$ for $1\le i_1<\ldots<i_l\le k$, $j>m$, therefore $H^j(X\setminus Y;\mathbb{Q})=0$ for $j\ge m+k$.\\

Now Lemma 3.3.1 implies that $h^{pq}(Y)=\delta_{pq}$ for $0\le p\le n$, $0\le q\le n$, $p+q\neq n$. So we get all $h^{pq}(Y)$ as soon as we know $e^p(Y)$, which can be computed using Theorem 2.6.4.\\

{\bf Examples:} Examples 3.2.4-3.2.6 can be treated much more easily by this method.\\

Another simple case is $\dim\,Y\le 1$. Note $Y$ is smooth because $Y$ is a $V$-manifold and that we have $h^{00}(Y)=b_0(Y)=$ number of connected components of $Y$, which settles already the case $\dim\,Y=0$. In the case $\dim\,Y=1$ it is therefore sufficient to compute the numbers $e^p(Y)$, using Theorem 2.3.5.\\

In the following examples we take $X=$ product of projective spaces in order to avoid the case $X=$ weighted projective space where Lemma 3.3.1 leads to a considerable simplification.\\

{\bf Example 3.3.2:}  The space $\mathbb{P}_2\times\mathbb{P}_1$ is a toric variety, it is a compactification of $\mathbb{C}^3$. Let $Y$ be the closure of $\{z\in\mathbb{C}^3\,|\,z_1^2+z_2^2=1, z_1-z_2z_3=0\}\simeq \{z\in\mathbb{C}^2\,|\,z_2^2(z_3^2+1)=1\}$. Solving the last equation for $z_2$ shows that $Y$ is connected, so using $e^p(Y)$ we get all Hodge numbers of $Y$.\\

Now there are other special arguments to simplify the computation.\\

For example, we can argue in the case of Example 3.3.2 as follows:\\
Let $Z$ be the closure of $\{z\in\mathbb{C}^2\,|\,z_1^2+z_2^2=1\}$ in $\mathbb{P}_2$. Then $Z$ is a rational curve, and $Y$ is the graph of the morphism $Z\to\mathbb{P}_1$: $(z_1,z_2)\mapsto\frac{z_1}{z_2}$. So $Y\simeq\mathbb{P}_1$.\\

{\bf Example 3.3.3:} The space $\mathbb{P}_3\times\mathbb{P}_1$ is a toric variety, it is a compactification of $\mathbb{C}^4$. Let $Y$ be the closure of $\{z\in\mathbb{C}^4\,|\,z_1^2+z_2^2+z_3^2=1, z_1-z_2z_4=0\}$. Then $Y$ is a smooth toric complete intersection.\\
Let $Z$ be the closure of $\{z\in\mathbb{C}^3\,|\,z_1^2+z_2^2+z_3^2=1\}$ in $\mathbb{P}_3$. The map $\{z\in\mathbb{C}^3\cap Z\,|\,z_2\neq 0\}\longrightarrow \mathbb{C}$: $(z_1,z_2,z_3)\longmapsto {z_1\over z_2}$ can be extended to a holomorphic map $Z\setminus\{(0,0,\pm1)\}\longrightarrow \mathbb{P}_1$. Then $Y$ is the closure of the graph of this mapping in $Z\times\mathbb{P}_1$. The fibre of $Y\mapsto Z$ over $(0,0,\pm1)$ is $\mathbb{P}_1$. So $h^{pq}(Y)=h^{pq}(Z)$ for $(p,q)\neq(1,1)$, $h^{11}(Y)=h^{11}(Z)+2$. For the Hodge numbers of $Z$, use Lemma 3.3.1.\\

{\bf Example 3.3.4:} Let $X:=\mathbb{P}_2\times\mathbb{P}_1$, $z_0,z_1,z_2$ and $w_0,w_1$ the homogeneous coordinates on $\mathbb{P}_2$ resp. $\mathbb{P}_1$, $g(z_0,z_1,z_2,w_0,w_1):=z_1w_0-z_2w_1$. Then $Y$ is smooth. We decompose $Y$:\\
$Y\cap\{w_0\neq 0,z_2\neq 0\}$ is the graph of the function $\mathbb{C}^2\to\mathbb{C}$: $(z_1,z_2)\mapsto w_1:=\frac{z_1}{z_2}$,\\
$Y\cap\{w_0\neq 0, z_2=0\}$ is isomorphic to $\mathbb{C}$,\\
$Y\cap\{w_0=0\}$ is isomorphic to $\mathbb{P}_1$.\\
Altogether we obtain that $Y$ has the same Hodge numbers as $\mathbb{P}_1\times\mathbb{P}_1$.\\
But $Y$ is not even homeomorphic to $\mathbb{P}_1\times\mathbb{P}_1$:\\
In fact $Y$ arises by blowing up the point $(1:0:0)$ in $\mathbb{P}_2$, see [Fi] p. 164.\\
At the same time, $Y$ is the Hirzebruch surface $\Sigma_1$. Now $\Sigma_0\simeq\mathbb{P}_1\times\mathbb{P}_1$, and by [Hi1] we have that $\Sigma_k$ is homeormorphic to $\Sigma_l$ if and only if $k-l$ is even.

\vspace{3cm}

\centerline{\bf References}
\vspace{1cm}

[A] M.Audin, The topology of torus actions on symplectic manifolds, Birkh\"auser: Basel 1991.\\

[Ba-C] V.V.Batyrev, D.A.Cox, {\sl On the Hodge structure of projective hypersurfaces in toric varieties}, Duke Math. J. 75, 293-338 (1994).\\

[C] D.A.Cox, {\sl The coordinate ring of a toric variety}, J. Algebraic Geom. {\bf 4}, 17-50 (1995).\\

[D] P.Deligne, {\sl Th\'eorie de Hodge II}, Publ. Math. IHES {\bf 40}, 5-58 (1971), {\sl III}, Publ. Math. IHES {\bf 44}, 5-77 (1974).\\

[Dz] T.Delzant, {\sl Hamiltoniens p\'eriodiques et image convexe de l'application moment}, Bull. Soc. Math. France {\bf 116} 315-339 (1988). \\

[D-K] V.I.Danilov, A.G.Khovanski\u{\i}, {\sl Newton polyhedra and an algorithm for computing Hodge-Deligne numbers}, Izv. Akad. Nauk SSSR Ser. Mat. {\bf 50}, no. 5 (1986), 925-945 = Math. USSR Izv. {\bf 29} (1987), 279-298.\\

[Fi] Complex analytic geometry. Lecture Notes in Mathematics {\bf 538}. Springer-Verlag: Berlin 1976.\\

[H1] H.A.Hamm, {\sl Rod $\chi_y$ kwasiodnorodnykh polnykh perese\v{c}eni\u{\i}}, Funkcional. Anal. i Prilo\v{z}en. {\bf 11} (1977), vyp. 1, 87-88 (Russian) = {\sl Genus $\chi_y$ of quasihomogeneous complete intersections}, Functional Anal. Appl. {\bf 11} (1977), 78-79.\\

[H2] H.A.Hamm, {\sl Hodge numbers for isolated singularities of nondegenerate complete intersections}, in: Singularities, Progress in Mathematics, Vol. 162, Birkh\"auser: Basel 1998, pp. 37-60.\\

[H3]  H.A.Hamm, {\sl Very good quotients of toric varieties}, in: Real and complex singularities, Chapman \& Hall/CRC: London 2000, pp. 61-75.\\

[Hi1] F. Hirzebruch: \"Uber eine Klasse von einfach zusammenh\"angenden komplexen Mannigfaltigkeiten. Math. Ann. {\bf 124}, 77-88 (1951).\\

[Hi2] F. Hirzebruch: Topological methods in algebraic geometry. Springer-Verlag: Berlin $^31966$.\\

[Kr] H. Kraft: Geometrische Methoden in der Invariantentheorie. Vieweg: Braunschweig 1984.\\

[K-K-M-S] G.Kempf, F.Knudsen, D.Mumford, B.Saint-Donat, {\sl Toroidal embeddings I}, Springer Lecture Notes {\bf 339} (1973).\\

[O] T.Oda, Convex bodies and algebraic geometry, Springer: Berlin 1988.\\

[Ok] M.Oka, {\sl A Lefschetz type theorem in a torus}. In: Singularity theory (Trieste 1991), pp. 574-593. World Sci. Publishing: River Edge, N.J. 1995.\\

[St] J.H.M. Steenbrink: Mixed Hodge structure on the vanishing cohomology. In: Real and Complex Singularities, Oslo 1976, pp. 525-563. Sijthoff\& Noordhoff: Alphen a/d Rijn 1977.\\

\end{document}